\documentclass[11pt,a4paper]{article}
\usepackage[english]{babel}  
\usepackage{geometry}
\usepackage{amsmath,amssymb,amsfonts}
\usepackage{mathrsfs}
\usepackage{color}
\usepackage{placeins}
\usepackage{multirow, makecell, array, enumerate}
\usepackage{colortbl}
\usepackage{bm} 
\usepackage{multicol}
\usepackage[width=.9\textwidth]{caption} 
\usepackage{dsfont}
\usepackage{graphicx}
\usepackage{subcaption}
\usepackage{algorithmicx}
\usepackage{algorithm}
\usepackage{algpseudocode}
\newtheorem{remark}{Remark}
\definecolor{gray}{rgb}{0.5, 0.5, 0.5}
\algdef{SE}[DOWHILE]{Do}{doWhile}{\algorithmicdo}[1]{\algorithmicwhile\ #1}

\definecolor{verde}{rgb}{0.2, 0.5, 0.3}
\newcommand\bl[1]{\textcolor{black}{#1}}

\title{A progressive reduced basis/empirical interpolation method for nonlinear parabolic problems\footnotemark[1]}            

\author{\scshape{Amina Benaceur\footnotemark[2] \footnotemark[3], 
Virgine Ehrlacher\footnotemark[2], Alexandre Ern\footnotemark[2] and 
S\'ebastien Meunier\footnotemark[3]}}

\date{}

\begin{document}
\maketitle 
\renewcommand{\thefootnote}{\fnsymbol{footnote}}
\footnotetext[1]{This work is partially supported by Electricit\'e De France (EDF) and a CIFRE PhD fellowship from ANRT}
\footnotetext[2]{University Paris-Est, CERMICS (ENPC), 77455 Marne la Vall{\'e}e Cedex 2 and INRIA Paris, 75589 Paris, France}
\footnotetext[3]{EDF Lab Les Renardi{\`e}res, 77250 Ecuelles Cedex, France}

\begin{abstract}
We investigate new developments of the combined 
Reduced-Basis and Empirical Interpolation Methods (RB-EIM) 
for parametrized nonlinear parabolic problems. 
In many situations, the cost of the EIM in the offline stage 
turns out to be prohibitive since 
a significant number of nonlinear
time-dependent problems need to be solved using the high-fidelity (or full-order) model. 
In the present work, we develop a new methodology, the Progressive  
RB-EIM (PREIM) method for nonlinear parabolic problems.
The purpose is to reduce the offline cost while
maintaining the accuracy of the RB approximation in the online stage.
The key idea is a progressive enrichment of both the EIM approximation 
and the RB space, in contrast to the standard approach where the EIM 
approximation and the RB space are built separately.
PREIM uses high-fidelity computations whenever available and RB computations
otherwise. Another key feature of each PREIM iteration is to select twice the
parameter in a greedy fashion, the second selection being made 
after computing the high-fidelity solution for the firstly selected
value of the parameter.
Numerical examples are presented on nonlinear heat transfer problems.
\end{abstract}

%

\section{Introduction}
The Reduced-Basis (RB) method devised in \cite{MMOPR:00,reliable} (see also
the recent textbooks~\cite{stamm16,manzoni16})
is a computationally 
effective approach to approximate parametrized Partial
Differential Equations (PDEs) encountered in many 
problems in science and engineering.
For instance, the RB method is often used in real-time simulations, where 
a problem needs to be solved very quickly under limited computational 
resources, or in multi-query simulations,
where a problem has to be solved repeatedly for a large
 number of parameter values.
Let $\mathcal P$ denote the parameter set. The RB method is
split into two stages: \textup{(i)} 
an offline stage where a certain number of so-called High-Fidelity (HF)
trajectories are computed for
a training subset of parameters 
$\mathcal P^\mathrm{tr}\subsetneq \mathcal P$(typically a finite 
element space based on a fine mesh); 
\textup{(ii)} an online stage for real-time or 
multi-query simulations where the parameter set $\mathcal P$
is explored more extensively. 
The output of the offline phase includes
an approximation space of small dimension spanned by the so-called RB functions.
The reduced space then replaces the much larger HF space in the online stage.
The crucial point for the computational efficiency of the overall procedure is that computations in the HF space are allowed only in the 
offline stage.

In the present work, we are interested in nonlinear parabolic problems
for which a RB method has been successfully developed in~\cite{grepl12,grepl07}.
A key ingredient to treat the nonlinearity so as to enable an online stage without HF 
computations is the Empirical Interpolation Method (EIM)~\cite{eim04, eim09}.
The EIM provides a separated approximation of the nonlinear (or non-affine) 
terms in the PDE. This approximation is built using a greedy algorithm as the 
sum of $M$ functions where the dependence on the space variable is separated
from the dependence on the parameter (and the time variable for parabolic problems). The integer $M$
is called the rank of the EIM and controls the 
accuracy of the approximation. 
Although the EIM is performed during the offline stage of the RB method,
its cost can become a critical issue since 
the EIM can require an important number of HF
computations for an accurate approximation of the nonlinearity.
The cost of the EIM typically scales with 
the size of the training set $\mathcal P^\mathrm{tr}$.

The goal of the present work is to overcome this issue.
To this purpose, we devise a new methodology, the Progressive RB-EIM (PREIM) method,
which aims at reducing the computational cost of the offline stage 
while maintaining
the accuracy of the RB approximation in the online stage. 
The key idea is a progressive enrichment of both the EIM approximation and the RB space, in contrast to the standard approach where the EIM 
approximation and the RB space are built separately. In PREIM, 
the number of HF computations is at most $M$,
and it is in general
much lower than $M$ in a time-dependent context where the greedy selection
of the pair $(\mu,k)$ to build the EIM approximation 
(where $\mu$ is the parameter and $k$ refers to the discrete time node) 
can lead to repeated values of $\mu$ for many different values of $k$. 
In other words, PREIM
can select multiple space fields within the same HF trajectory to 
build the EIM space functions.  
In this context, only a modest number of HF trajectories needs to be computed,
yielding significant computational savings with respect to the standard offline
stage. \bl{PREIM is driven by convergence criteria on the quality of both the EIM and the RB approximation, as in the standard RB-EIM procedure. Moreover, as the PREIM iteration progresses, more and more HF functions are used to approximate the nonlinearity, thus attaining the same quality as the standard approach at convergence. The benefit of PREIM is to deliver accurate EIM and RB approximations at moderate offline costs. This is possible in situations where the computation of HF trajectories dominates the cost of the progressive construction of the EIM and the RB, which is observed to be the case in our numerical experiments. Moreover, PREIM is expected to bring computational savings whenever the nonlinearity can be represented by a separated approximation of relatively modest rank since otherwise PREIM will need to compute a number of HF trajectories comparable to that of the standard procedure.}

The idea of a progressive enrichment of both the EIM approximation 
and the RB space has been recently proposed
in~\cite{ser15} for stationary nonlinear PDEs, where it is called
Simultaneous EIM/RB (SER). Thus, PREIM extends this idea to time-dependent
PDEs. In addition, there is an important difference in the greedy algorithms 
between SER and PREIM. 
Whereas SER uses only RB computations, PREIM uses HF
 computations whenever available, both for the greedy selection 
of the parameters and the time nodes, as well as for the 
space-dependent functions 
in the EIM approximation. These aspects are particularly relevant
since they improve the accuracy of the
EIM approximation. This is illustrated in our numerical
experiments on nonlinear parabolic PDEs. 
\bl{The progressive construction of the EIM and the RB has been
recently addressed within the Empirical Interpolation Operator Method
in~\cite{drohmann}. Therein, the enrichment criterion is common to both the EIM and the RB, and the snapshot maximizing an a posteriori error estimator is selected to enrich both bases. Instead, PREIM has dedicated criteria for the quality
of the EIM approximation and for the RB approximation. Furthermore, PREIM systematically exploits the knowledge of the HF trajectories whenever available, and an update step is performed so as to assess the current parameter selection.
We also mention the Proper Orthogonal Empirical Interpolation Method from~\cite{poim} where a progressive construction of the EIM approximation is devised using POD-based approximations of the HF trajectories. 
}
 
The paper is organized as follows.
In Section~\ref{mod_pb}, we introduce the model problem. 
In Section~\ref{rb_strat}, we briefly recall the main ideas of the 
nonlinear RB method devised in~\cite{grepl12,grepl07}, and 
in Section~\ref{off}, we briefly recall the EIM procedure in the standard
offline stage as devised in~\cite{eim04,eim09}. 
The reader familiar with the material can jump directly to
Section~\ref{sec:PREIM} where PREIM is introduced and discussed. 
Section~\ref{num_res} presents numerical results
illustrating the performance of PREIM on nonlinear parabolic
problems related to heat transfer. Finally,
Section~\ref{sec:conc} draws some conclusions and outlines some
perspectives.

\section{Model problem}
\label{mod_pb}
In this section, we present a prototypical example of a nonlinear parabolic PDE.
The methodology we propose is illustrated on this model problem but can be extended to other types of parabolic equations.
We consider a spatial domain (open, bounded, 
connected subset) $\Omega \subset \mathbb{R}^d$, $d \geq 1$, 
with a Lipschitz boundary, a finite time interval $I=[0,T]$, with $T>0$, 
and a parameter set $\mathcal{P} \subset \mathbb{R}^{p}$, $p \geq 1$,
whose elements are generically denoted by $\mu \in \mathcal{P}$.
Our goal is to solve the following nonlinear parabolic PDE for many values 
of the parameter \bl{$\mu\in\mathcal P$}: find 
$u_\mu:I\times \Omega\rightarrow \mathbb R $ such that
  \begin{align}\label{HTE}
  \left\{
  \begin{alignedat}{2}
    \frac{\partial u_\mu}{\partial t}
   - \nabla \cdot \big((\kappa_0 + \Gamma(\mu,u_\mu)) \nabla u_\mu\big)& 
   = f, &&\mathrm{in}\ I\times\Omega , \\
  -\big(\kappa_0 + \Gamma(\mu,u_\mu)\big) \frac{\partial u_\mu}{\partial n}&
  = \phi_e,\ \ \ \ &&\mathrm{on}\  I\times\partial\Omega , \\
  u_\mu(t=0,\cdot)&= u_0(\cdot),\quad && \mathrm{in}\  \Omega,
  \end{alignedat}
  \right.
  \end{align}
where $\kappa_0 >0$ is a fixed positive real number, 
$\Gamma:\mathcal{P}\times \mathbb R \rightarrow \mathbb{R}$ is 
a given nonlinear function,
$f:I\times\Omega\rightarrow\mathbb R$ is the source term,
$\phi_e:I\times\partial\Omega\rightarrow\mathbb R$ is the time-dependent \bl{
Neumann boundary condition} on $\partial\Omega$,
and $u_0:\Omega\rightarrow \mathbb R$ is the initial 
condition. For simplicity, we assume without loss of generality
that $f$, $\phi_e$, and $u_0$
are parameter-independent. 
We assume that $f\in L^2(I;L^2(\Omega))$ and $\phi_e\in L^2(I;L^2(\partial\Omega))$ (this means 
that $f(t)\in L^2(\Omega)$ and $\phi_e(t)\in L^2(\partial\Omega)$ for (almost every) $t\in I$), 
and we also assume that $u_0\in H^1(\Omega)$.
We make the standard uniform ellipticity assumption
$\beta_1 \leq \kappa_0 + \Gamma(\mu,z) \leq \beta_2$
with $0<\beta_1<\beta_2<\infty$, for all $(\mu,z) \in \mathcal P
\times \mathbb R$. We do not specify a functional setting for the 
nonlinear parabolic PDE~\eqref{HTE}; with the above assumptions, it is 
reasonable to look for a weak solution $u\in L^2(I;Y)$ with $Y = H^{1}(\Omega)$.

\begin{remark}[Initial condition]
For parabolic PDEs, the initial condition is often taken to be in a
larger space, e.g., $u_0\in L^2(\Omega)$. Our assumption that $u_0\in Y$
is motivated by the RB method where basis functions in $Y$ are sought
as solution snapshots in time and for certain parameter values.
In this context, we want to include the possibility to select the initial condition as a RB function. 
\end{remark}

\begin{remark}[Heat transfer] \label{rem:heat}
One important application we have in mind for~\eqref{HTE} is heat transfer 
problems. In this context, the PDE can take the slightly more general form
\[
\alpha(u_\mu) \frac{\partial u_\mu}{\partial t}
- \nabla \cdot \big((\kappa_0 + \Gamma(\mu,u_\mu)) \nabla u_\mu\big) 
= f,\quad \mathrm{in}\ I\times\Omega,
\]
where $\alpha(u_\mu)$ stands for
the mass density times the heat capacity. Moreover, the quantity
$(\kappa_0 + \Gamma(\mu,u_\mu))$ represents the thermal 
conductivity. Note also that $\phi_e>0$ means that the system is heated.
\end{remark}

In practice, one way to solve~\eqref{HTE} is to use a $Y$-conforming
Finite Element Method~\cite{ern_guermond} to discretize in space and 
a time-marching scheme to discretize in time. The Finite Element 
Method is based on 
a finite element subspace $X \varsubsetneq Y$
defined on a discrete nodal subset $\Omega^\mathrm{tr}\varsubsetneq\Omega$, 
where $\mathrm{Card}(\Omega^\mathrm{tr}) = \mathcal N$.
To discretize in time, we consider an integer $K\geq1$, we let
$ 0=t^0<\cdots<t^K=T$ be $(K+1)$ distinct
time nodes over $I$, and we set $ \mathbb K^\mathrm{tr} = \{1,\cdots,K\}$, 
$ \overline {\mathbb K}^\mathrm{tr} = \{0\}\cup\mathbb K^\mathrm{tr}$, 
$I^\mathrm{tr} = \{t^k\}_{k\in \overline{\mathbb K}^\mathrm{tr}}$,
and $ \Delta t^k = t^k - t^{k-1}$ for all $k \in \mathbb K^\mathrm{tr}$.
As is customary with the RB method, we assume henceforth that the mesh-size
and the time-steps are small enough so that the above space-time
discretization method delivers HF approximate trajectories within the
desired level of accuracy. These trajectories, 
which then replace the exact trajectories solving~\eqref{HTE}, 
are still denoted $u_\mu$ for all $\mu\in\mathcal P$.
Henceforth, we use the convention that the superscript $k$
always indicates a time index; thus, we write
$u^k_\mu(\cdot) = u_\mu(t^k,\cdot)\in X$, $f^k(\cdot) = f(t^k,\cdot)\in L^2(\Omega)$, and $\phi_e^k(\cdot) = \phi_e(t^k, \cdot)\in L^2(\partial\Omega)$. 
Applying a semi-implicit Euler scheme, 
our goal is, given $u^0_\mu = u_0\in X$, to find
$(u^k_{\mu})_{k\in\mathbb K^\mathrm{tr}} \in X^K$ such that,
for all $k \in \mathbb K^\mathrm{tr}$,
\begin{equation}\label{tnl1}
\begin{alignedat}{2}
\forall v \in X, \enskip
m( u^k_\mu,v)
 + \Delta t^k  a_0\big(u^k_\mu, v\big)
 + \Delta t^k n_\Gamma\big(\mu,u^{k-1}_\mu, v\big)
= m( u^{k-1}_\mu, v) 
 + \Delta t^k l^k(v), 
\end{alignedat}
\end{equation}
with the bilinear forms 
$m: Y\times Y\rightarrow \mathbb R$, $a_0: Y\times Y\rightarrow \mathbb R$ 
and the linear forms $l^k: Y\rightarrow \mathbb R$ such that
\begin{equation}\label{ma0}
m(v,w) =\int_\Omega vw,
\qquad
a_0(v,w) =\kappa_0\int_\Omega \nabla v\cdot\nabla w,
\qquad
l^k(v) = \int_\Omega f^k v + \int_{\partial\Omega}\phi_e^k v,
\end{equation}
and the nonlinear form
$n_\Gamma:\mathcal P \times  Y\times Y\rightarrow \mathbb R$ such that
\begin{equation}
n_\Gamma(\mu,v,w) = \int_\Omega \Gamma(\mu,v) \nabla v\cdot \nabla w,
\end{equation}
for all $\mu \in \mathcal P$ and all $v,w \in Y$. 
\bl{In~\eqref{tnl1}, the nonlinearity is treated
explicitly, whereas the diffusive term is treated implicitly.
This choice avoids dealing
with a nonlinear solver at each time-step. 
The computation of derivatives of discrete operators within Newton's method
is addressed, e.g., in \cite{drohmann}.}

\section{The Reduced-Basis method}\label{rb_strat}
In this section, we briefly recall 
the Reduced-Basis (RB) method for the nonlinear 
problem~\eqref{tnl1}~\cite{grepl07,grepl12}. 
Let $\hat X_N \subset X$ be a so-called reduced subspace such that
$ N=\mathrm{dim}(\hat X_N) \ll \mathrm{dim}(X) = \mathcal N$. Let 
$(\theta_n)_{1\leq n \leq N}$ be a $Y$-orthonormal basis of $\hat X_N$. 
For all $\mu\in\mathcal P$ and 
$k \in \overline{\mathbb K}^\mathrm{tr}$, the RB solution  
$\hat{u}^k_{\mu}\in \hat X_N$ that approximates the HF solution
$u^k_\mu\in X$ is decomposed as
\begin{equation}\label{rbsol}
\hat{u}^k_{\mu}=\sum_{n=1}^N \hat{u}^k_{\mu,n}
\theta_n,
\end{equation} 
with real numbers $\hat u ^k_{\mu,n}$ for all $n\in\{1,\ldots,N\}$.
Let us introduce the component vector 
$\mathbf{\hat{u}}^k_{\mu}:=(\hat u ^k_{\mu,n})_{1\leq n\leq N} \in \mathbb R^N$, for all $\mu \in \mathcal P$ and $k \in \overline{\mathbb K}^\mathrm{tr}$.
Let $\hat u^0$ be the $Y$-orthogonal \bl{projection} of the initial condition $u_0\in X$ onto $\hat X_N$ with associated component vector $\mathbf {\hat{u}}^0\in\mathbb R^N$. Replacing $u^k_\mu\in X$
in the weak form~\eqref{tnl1} by the approximation $\hat{u}^k_{\mu}\in \hat X_N$ with associated component vector $\mathbf{\hat{u}}^k_{\mu}\in \mathbb R^N$, and using the test functions $(\theta_p)_{1\leq p \leq N}$,
we obtain the following problem written in algebraic form: 
Given $\mathbf{\hat{u}}^0_{\mu} = \mathbf {\hat{u}}^0\in\mathbb R^N$, 
find $(\mathbf{\hat{u}}^k_{\mu})_{k\in {\mathbb K}^\mathrm{tr}}\in (\mathbb R^N)^K$ such that, for all $k\in {\mathbb K}^\mathrm{tr}$,
\begin{align}\label{rbmat0}
(\mathbf M  + \Delta t^k {\mathbf A_0}) \mathbf{\hat{u}}^{k}_{\mu} 
 = \Delta t^k  \mathbf f^k 
 + \mathbf M \mathbf{\hat{u}}^{k-1}_{\mu}
 -\Delta t^k  \mathbf g(\mathbf{\hat{u}}^{k-1}_{\mu}),
\end{align}
with the matrices $\mathbf M, \mathbf A_0 \in \mathbb R^{N\times N}$ and
the vectors $\mathbf f^k\in\mathbb R^N$ such that
\begin{align}\label{matrices}
\mathbf M = \Big(m(\theta_n,\theta_p)\Big)_{1\leq p,n \leq N }, \qquad
\mathbf A_0 = \Big(a_0(\theta_n,\theta_p)\Big)_{1\leq p,n \leq N }, \qquad
\mathbf f^k = \big(l^k(\theta_p)\big)_{1\leq p \leq N },
\end{align}
and the vector ${\mathbf g}(\mathbf{\hat{u}}^{k-1}_{\mu}) \in \mathbb R^N$ such that
\begin{equation} \label{eq:G_bad}
{\mathbf g}(\mathbf{\hat{u}}^{k-1}_{\mu})
= \bigg( \sum_{n=1}^N\hat u^{k-1}_{\mu,n} \int_\Omega\Gamma\bigg(\mu,\sum_{n'=1}^N \hat u^{k-1}_{\mu,n'}\theta_{n'}\bigg)
\nabla \theta_n\cdot\nabla\theta_p \bigg)_{1\leq p \leq N }.
\end{equation}
The difficulty is that 
the computation of ${\mathbf g}(\mathbf{\hat{u}}^{k-1}_{\mu})$ 
requires a parameter-dependent reconstruction using the RB functions so as to compute the integral over $\Omega$.
To avoid this, we need to build an approximation $\gamma_M$ of the nonlinear function
$\gamma:\mathcal{P}\times \overline{\mathbb K}^\mathrm{tr}\times \Omega
\rightarrow \mathbb{R}$ such that
\begin{equation}\label{sepa}
\gamma(\mu,k,x) := \Gamma(\mu,u_\mu^k(x)),
\end{equation}
in such a way that the dependence on
$x$ is separated from the dependence on $(\mu,k)$.
More precisely, for some integer $M>0$, we are looking for an (accurate) approximation 
$\gamma_M:\mathcal{P}\times\overline{\mathbb K}^\mathrm{tr}\times\Omega\rightarrow \mathbb{R}$ 
of $\gamma$ under the separated form
\begin{equation}\label{eim}
\gamma_M(\mu,k,x)  := \sum_{j=1}^M \varphi_{\mu,j}^{k}  q_j(x),
\end{equation}
where $M$ is called the rank of the approximation and $\varphi^k_{\mu,j}$ are 
real numbers
\bl{
that we find by interpolation over a set of $M$ points $\{x_1, \ldots, x_M\}$ in $\Omega^\mathrm{tr}$ by requiring that
\begin{equation} \label{eq:interpolation}
\gamma_M(\mu, k, x_i ) = \gamma(\mu, k, x_i) = \Gamma(\mu,u^{k}_\mu(x_i)),\quad \forall i \in \{1,\cdots,M\}.
\end{equation}
The interpolation property~\eqref{eq:interpolation} is achieved
by setting
\begin{equation}
\varphi_{\mu,j}^{k} = (\mathbf B^{-1} \bm{\gamma}_{\mu}^{k})_j, \; \forall j \in \{1,\cdots,M\},
\quad \text{where} \quad \bm{\gamma}_{\mu}^{k} :=\big( \Gamma(\mu,u^{k}_\mu(x_i))
\big)_{1\leq i \leq M}\in \mathbb R^M,
\end{equation}
and $\mathbf B=(q_j(x_i))_{1\le i,j\le M} \in \mathbb R^{M\times M}$ must be an invertible matrix.
Therefore, \eqref{eim} can be rewritten as follows:
\begin{equation} \label{eq:eim_bis}
\gamma_M(\mu,k,x) = \sum_{j=1}^M (\mathbf B^{-1} \bm{\gamma}_{\mu}^{k})_j q_j(x).
\end{equation}
The points $(x_i)_{1\leq i \leq M}$ in $\Omega^\mathrm{tr}$ and
the functions $(q_j)_{1\leq j \leq M}$ defined on $\Omega$
are determined by the EIM algorithm \cite{eim04} which is further described in Section~\ref{off} below.}

\bl{Let us now describe how we can use the EIM approximation~\eqref{eq:eim_bis} to allow for an offline/online decomposition of the computation of the vector ${\mathbf g}(\mathbf{\hat{u}}^{k-1}_{\mu})$ defined in~\eqref{eq:G_bad}.
Under the (reasonable) assumptions $\hat{u}^{k}_{\mu} \approx u_\mu^{k}$ and
$\Gamma(\mu,\hat{u}^{k}_{\mu}(x)) \approx
\Gamma\big(\mu,u_\mu^{k}(x)\big)$, we obtain
\begin{align}
\Gamma(\mu,\hat{u}^{k}_{\mu}(x)) &\approx
\Gamma\big(\mu,u_\mu^{k}(x)\big) = \gamma(\mu,{k},x)
\approx \gamma_M(\mu,{k},x) 
=\sum_{j=1}^M (\mathbf B^{-1} \bm{\gamma}_{\mu}^{k})_j q_j(x)
\approx \sum_{j=1}^M (\mathbf B^{-1} \bm{\hat \gamma}_{\mu}^{k})_j q_j(x)
\end{align}
with the vector $\bm{\hat \gamma}_{\mu}^{k}:=( \Gamma(\mu,\hat{u}^{k}_\mu(x_i)))_{1\leq i \leq M}\in\mathbb{R}^M$.}
The problem~\eqref{rbmat0} then becomes:
Given $\mathbf{\hat{u}}^0_{\mu} = \mathbf {\hat{u}}^0\in \mathbb R^N$, 
find $(\mathbf{\hat{u}}^k_{\mu})_{k\in\mathbb K^\mathrm{tr}}\in(\mathbb R^N)^K$ such that, for all $k\in \mathbb K^\mathrm{tr}$,
\begin{align}\label{rbmat}
(\mathbf M  + \Delta t^k {\mathbf A_0}) \mathbf{\hat{u}}^{k}_{\mu} 
 = \Delta t^k  \mathbf f^k 
 + \big(\mathbf M 
 -\Delta t^k  \mathbf D^{k-1}_\mu\big)\mathbf{\hat{u}}^{k-1}_{\mu},
\end{align} 
with the matrix $\mathbf D^{k-1}_\mu \in \mathbb R^{N\times N}$ such that
\begin{equation}\label{matrices2}
\mathbf D^{k-1}_\mu
= \sum_{j=1}^M\big(\mathbf B^{-1}\bm{\hat \gamma}_\mu^{k-1}\big)_j \mathbf C^j \quad \mbox{where}\quad 
\mathbf C^j = \bigg(\int_\Omega q_j\nabla\theta_n\cdot \nabla \theta_p \bigg)_{1\leq p,n \leq N }\in \mathbb{R}^{N\times N}, \; \forall 1\leq j \leq M.
\end{equation}

The overall computational procedure can now be split into
two stages: 
\begin{enumerate}[(i)]
\item An offline stage where one precomputes on the one hand 
the RB functions $(\theta_n)_{1\leq n \leq N}$ leading to the 
vectors $\mathbf{\hat{u}}^0\in \mathbb R^N$,
$(\mathbf f^k)_{k\in \mathbb K^\mathrm{tr}}\in (\mathbb R^N)^K$ and the matrices
$\mathbf M,\mathbf A_0\in \mathbb R^{N\times N}$, and on the other hand
the EIM points $(x_i)_{1\leq i \leq M}$ and
the functions $(q_j)_{1\leq j \leq M}$ leading 
to the matrices 
$ \mathbf B\in \mathbb R^{M\times M}$ and $\mathbf C^j\in \mathbb R^{N\times N}$,
for all $j\in\{1,\ldots,M\}$. The offline stage is discussed in more
detail in Section~\ref{off}.
\item An online stage to be performed each time one wishes 
to compute a new trajectory for a 
parameter $\mu \in \mathcal P$. All what remains to be performed  
is to compute the vector $\widehat{\bm{\gamma}}_\mu^{k-1}\in \mathbb R^M$ and the matrix
$\mathbf D^{k-1}_\mu\in \mathbb R^{N\times N}$ and to 
solve the $N$-dimensional linear problem~\eqref{rbmat} for all $k\in \mathbb K^\mathrm{tr}$. 
The online stage is summarized in Algorithm~\ref{online}.
\end{enumerate}
\begin{algorithm}[h]
\caption{Online stage\label{online}}
\begin{algorithmic}[1]
\vspace{0.2cm}
\Statex {\scshape{\bfseries \underline  {Input :}}} 
\bl{$\mu $,} $(\theta_n)_{1\leq n\leq N}$, $\mathbf{\hat{u}}^0$, 
 $(\mathbf f^k)_{k\in{\mathbb K}^\mathrm{tr}}$, $\mathbf M$, $\mathbf A_0$,
 $(x_i)_{1\leq i \leq M}$, $(q_j)_{1\leq j \leq M}$, and $(\mathbf C^j)_{1\leq j\leq M}$
\State Set $k=1$ and
$\mathbf{\hat{u}}^0_{\mu} = \mathbf{\hat{u}}^0$
\While {$k\in \mathbb K^\mathrm{tr}$ }
\State Compute 
$\mathbf D_{\mu}^{k-1}$ using~\eqref{matrices2} and $\mathbf{\hat{u}}^{k-1}_\mu$
\State Solve the reduced system \eqref{rbmat} 
to obtain $\mathbf{\hat{u}}^k_\mu$
\State Set $k = k+1$
\EndWhile
\Statex {\scshape{\bfseries \underline  {Output :}}}
$(\mathbf{\hat{u}_{\mu}}^{k})_{k\in\mathbb K^\mathrm{tr}}$
\vspace{0.2cm}
\end{algorithmic}
\end{algorithm}

\section{The standard offline stage}\label{off}

There are two tasks to be performed during 
the offline stage: 
\bl{
\begin{enumerate}[(T$_1$)]
\item Build the rank-$M$ EIM
approximation~\eqref{eim}  
of the nonlinear function $\gamma$ defined
by~\eqref{sepa};
\item Explore the solution manifold in order to construct 
a linear subspace $\hat X_N \subset X$ of dimension $N$. 
\end{enumerate}}
In the standard offline stage, these two tasks are performed independently.

\begin{algorithm}[ht]
\caption{Standard EIM\label{offline_eim}}
\begin{algorithmic}[1]
\vspace{0.2cm}
\Statex {\scshape{\bfseries \underline  {Input :}}}  
$\mathcal{P}^\mathrm{tr}$,
$\overline{\mathbb K}^\mathrm{tr}$, $\Omega^\mathrm{tr}$,
and $\epsilon_\textsc{eim}>0$
\State Compute $\mathcal S
= (u^k_\mu)_{(\mu,k)\in\mathcal{P}^\mathrm{tr}\times \overline{\mathbb K}^\mathrm{tr}}$
\label{off_comp}\hfill{\small $P$ \verb|HF trajectories|} 
\State Set $m=1$ and $\gamma_0\equiv0$
\State Search $(\mu_m, k_m) \in
 \underset{(\mu,k)\in \mathcal{P}^\mathrm{tr} 
\times \overline{\mathbb K}^\mathrm{tr}} {\operatorname{argmax}}
\|\Gamma(\mu, u^k_\mu(\cdot))
-\gamma_{m-1}(\mu,k,\cdot)\|_{\ell^\infty(\Omega^\mathrm{tr})}$
\State Set $r_m(\cdot)
:= \Gamma(\mu_m,u^{k_m}_{\mu_m}(\cdot)) 
-\gamma_{m-1}(\mu_m,{k_m},\cdot)$ and $x_m 
\in \underset{x \in \Omega^\mathrm{tr}}{\mathrm{argmax}}\ |r_m(x)|$
\While{($|r_m(x_m)|>\epsilon_\textsc{eim})$}
\State Set $q_m := {r_m}/{r_m(x_m)}$ and compute $(\mathbf B_{mi})_{1\leq i \leq m}$ by setting $\mathbf B_{mi} := (q_i(x_m))$ \label{8}
\State Set $m = m+1$
\State Search $(\mu_m, k_m) \in
 \underset{(\mu,k)\in \mathcal{P}^\mathrm{tr} 
\times \overline{\mathbb K}^\mathrm{tr}} {\operatorname{argmax}}
\|\Gamma(\mu, u^k_\mu(\cdot))
-\gamma_{m-1}(\mu,k,\cdot)\|_{\ell^\infty(\Omega^\mathrm{tr})}$
\State Set $r_m(\cdot) := \Gamma(\mu_m,u^{k_m}_{\mu_m}(\cdot)) 
-\gamma_{m-1}(\mu_m,{k_m},\cdot)$ and $x_m 
\in \underset{x \in \Omega^\mathrm{tr}}{\mathrm{argmax}}\ |r_m(x)|$\label{rm_xm_eim}
\EndWhile
\State Set $M:=m-1$
\Statex {\scshape{\bfseries \underline {Output :}}} 
$(x_i)_{1\le i\le M}$ and $(q_j)_{1\le j\le M}$
\vspace{0.2cm}
\end{algorithmic}
\end{algorithm}
Let us first discuss \bl{Task (T$_1$)}, i.e., the construction of the 
rank-$M$ EIM approximation. Recall from Section~\ref{rb_strat} that 
the goal is to find the interpolation points $(x_i)_{1\le i\le M}$ in $\Omega^\mathrm{tr} \varsubsetneq \Omega$ and the functions $(q_j)_{1\le j\le M}$ 
where $q_m:\Omega\rightarrow \mathbb R$.
The construction of the EIM approximation additionally uses
a training set $\mathcal P^\mathrm{tr} \subset \mathcal P$ for the 
parameter values; in what follows, we denote by $P$ the cardinality of 
$\mathcal P^\mathrm{tr}$.
For a real-valued function $v$ defined on $\Omega^\mathrm{tr}$, we define 
$\|v\|_{\ell^\infty(\Omega^\mathrm{tr})}:=
\mathrm{max}\bl{_{x\in \Omega^\mathrm{tr}}|v(x)|}$.
Given an iteration counter $m\geq1$ and a function $\gamma_{m-1}$ defined on $\mathcal P^\mathrm{tr}\times \overline{\mathbb K}^{\mathrm{tr}}\times \Omega$, with the convention that $\gamma_0\equiv0$, 
an EIM iteration consists of the following steps. First, one defines 
$(\mu_{m}, k_{m})\in\mathcal P^{\mathrm{tr}}\times\overline{\mathbb K}^\mathrm{tr}$ by
\begin{equation}\label{std_greedy}
(\mu_{m}, k_{m}) \in
\underset{(\mu,k) \in\mathcal P^{\mathrm{tr}}\times \overline{\mathbb K}^\mathrm{tr}}
{\mathrm{argmax}}\ \| \Gamma(\mu,u^k_\mu(\cdot)) 
- \gamma_{m-1}(\mu,k,\cdot)\|_{\ell^\infty(\Omega^\mathrm{tr})} ,
\end{equation}
where we notice the use of the HF trajectories for all values 
of the parameter $\mu$ in the 
training set $\mathcal P^\mathrm{tr}$.
Once $(\mu_{m}, k_{m})$ has been determined, one sets 
\begin{equation}\label{resid}
r_{m}(\cdot) := \Gamma(\mu_{m},u_{\mu_{m}}^{k_{m}}(\cdot)) 
 - \gamma_{m-1}(\mu_{m},{k_{m}},\cdot), \qquad
x_{m} \in  \underset{x\in\Omega^\mathrm{tr}}{\mathrm{argmax}}\ |r_m(x)|, 
\end{equation}
and one checks whether $|r_m(x_m)|>\epsilon_\textsc{eim}$
for some user-defined positive threshold $\epsilon_\textsc{eim}$.
If this is the case, one sets
\begin{equation}
q_{m}(\cdot) := \frac{r_{m}(\cdot)}{r_{m}(x_{m})},
\end{equation}
and one computes the new row of the matrix $\mathbf B$ by setting
$\mathbf B_{mi} := (q_i(x_m))$, for all $1\leq i \leq m$.
The standard EIM procedure is presented in Algorithm~\ref{offline_eim}.

Let us now briefly discuss \bl{Task (T$_2$)} above, i.e., the construction of 
a set of RB functions with cardinality $N$. 
First, as usual
in RB methods, the solution manifold 
is explored by considering a training set
for the parameter values; for simplicity, we consider the same
training set $\mathcal P^{\mathrm{tr}}$ as for the EIM approximation. 
This way, one only explores
the collection of points $\{u^k_{\mu}\}_{(\mu,k)\in\mathcal P^{\mathrm{tr}}\times \overline{\mathbb K}^\mathrm{tr}}$ in the solution manifold. For this exploration to be informative, the training set $\mathcal P^{\mathrm{tr}}$ has to be chosen large enough. The exploration can be driven by means of an a
posteriori error estimator (see, e.g.,~\cite{rovas06}) which allows one to
evaluate only $N$ HF trajectories.
However, in the present setting where HF trajectories are to be computed
for all the parameters in $\mathcal P^{\mathrm{tr}}$ 
when constructing the EIM approximation,
it is natural to exploit these computations by means 
of a Proper Orthogonal Decomposition (POD)~\cite{Hinze_05,podref} 
to define the RB. This technique is often considered
in the literature to build the RB in a time-dependent 
setting, see, e.g., \cite{haasdonk13,stamm16,manzoni16}.
\bl{
In practice, a POD of the whole collection of snapshots may be computationally 
demanding (or even unfeasible) when a very large number of functions is considered. 
Thus, we adopt a 
POD-based progressive construction of the reduced basis
in the spirit of the POD-greedy algorithm from~\cite{haasdonk13}. Therein, one additional RB function is picked at a time, whereas here we can pick more than one function. 
The progressive  construction of the RB is presented in 
Algorithm~\ref{offline_rb} where we have chosen an enumeration of the parameters in $\mathcal{P}^\mathrm{tr}$ from $1$ to $P$. 
The initialization of Algorithm~\ref{offline_rb} is made by computing $(\theta_n)_{1\leq n \leq N^1}=\mathrm{POD}(\mathcal S_1,\epsilon_\textsc{pod})$ 
for the trajectory $\mathcal S_1$ associated with the parameter $\mu_1$. That is, we select the first $N^1$ POD modes out of the set $\mathcal S_1$ with error threshold $\epsilon_\textsc{pod}$ (for completeness, this procedure
is briefly outlined in Appendix~\ref{pod_app}). The next steps of the algorithm are performed in an iterative fashion. For each new trajectory, we first subtract its projection on the current RB, and then perform a POD on the projection and merge the result with the current RB. This specific part of the procedure, called UPDATE\_RB, is presented in Algorithm~\ref{progressive_pod}; this part of the procedure is presented separately since it will be re-used later on.}

\begin{algorithm}[h]
\caption{\bl{Progressive RB}\label{offline_rb}}
\begin{algorithmic}[1]
\vspace{0.2cm}
\Statex {\scshape{\bfseries \underline  {Input :}}}  
$\mathcal{P}^\mathrm{tr}$,
$\overline{\mathbb K}^\mathrm{tr}$, and $\epsilon_\textsc{pod}>0$
\State Compute $\big(\mathcal S_p\big)_{1\leq p \leq P}
= \big((u^k_{\mu_p})_{k\in\overline{\mathbb K}^\mathrm{tr}}\big)_{1\leq p \leq P}$
\label{line1}\hfill {\small $P$ \verb|HF trajectories|}
\State Compute $(\theta_n)_{1\leq n \leq N^1}=\mathrm{POD}(\mathcal S_1,\epsilon_\textsc{pod})$\label{line:ini}
\State Set $p=1$
\While {$p< P$}
\State Set $p=p+1$
\State Compute $(\theta_n)_{1\leq n \leq N^p}= \mathrm{UPDATE\_RB}
\big((\theta_n)_{1\leq n\leq N^{p-1}}, \mathcal S_p, \epsilon_\textsc{pod}\big)$\label{5_pod}
\EndWhile
\State Set $N:=N^P$
\State Compute $\hat{\mathbf u}^0$, $(\mathbf f^k)_{k\in\mathbb K^\mathrm{tr}}$, 
$\mathbf M$, and $\mathbf A_0$
\State Compute the matrices $(\mathbf C^j)_{1\leq j \leq M}$
\Statex {\scshape{\bfseries \underline {Output :}}} $(\theta_n)_{1\leq n \leq N}$,
$\hat{\mathbf u}^0$,  $(\mathbf f^k)_{k\in\mathbb K^\mathrm{tr}}$, $\mathbf M$, $\mathbf A_0$, 
and $(\mathbf C^j)_{1\leq j \leq M}$
\vspace{0.2cm}
\end{algorithmic}
\end{algorithm}

\begin{algorithm}[htb]
\caption{\bl{UPDATE\_RB}\label{progressive_pod}}
\begin{algorithmic}[1]
\vspace{0.2cm}
\Statex {\scshape{\bfseries \underline {Input :}}} 
$\Theta=(\theta_n)_{1\le n\le N}$, $\mathcal{S}$, and $\epsilon_\textsc{pod}>0$
\If{$\mathcal{S}=\emptyset$}
\State $\Theta$ remains unchanged
\Else
\State Define $\tilde{\mathcal S}
:=(u - \Pi_{\mathrm{span}(\Theta)} u)_{u\in \mathcal{S}}$
\State Set $\Xi := \mathrm{POD}(\tilde{\mathcal S},\epsilon_\textsc{pod})$
\If{$\Xi =\emptyset$}
\State $\Theta$ remains unchanged
\Else
\State Set $\Theta := \Theta\cup \Xi$
\EndIf 
\EndIf
\Statex {\scshape{\bfseries \underline {Output :}}} $\Theta$
\vspace{0.2cm}
\end{algorithmic}
\end{algorithm}

\begin{remark}[Threshold $\epsilon_\textsc{pod}$]
\bl{For the initialization (line~\ref{line:ini} of Algorithm~\ref{offline_rb}),
one can use a relative error threshold for $\epsilon_\textsc{pod}$ (for instance, $\epsilon_\textsc{pod}=1\%$). 
Instead, for the iterative loop (line~\ref{5_pod} of Algorithm~\ref{offline_rb}), the threshold $\epsilon_\textsc{pod}$ can 
be set to the greatest singular value that has been truncated at the 
initialization step.}
\end{remark}
\begin{remark}[Order of EIM and RB]
Algorithms~\ref{offline_eim} 
and~\ref{offline_rb} can be performed in whatever order. 
If Algorithm~\ref{offline_rb} is performed first, 
the computation of the matrices $(\mathbf C^j)_{1\leq j \leq M}$ is postponed 
to the end of Algorithm~\ref{offline_eim}. 
Moreover, the HF trajectories 
$(u^k_\mu)_{(\mu,k)\in\mathcal{P}^\mathrm{tr}\times
\overline{\mathbb K}^\mathrm{tr}}$ appearing in both algorithms 
are computed only once.
\end{remark}

\section{The Progressive RB-EIM method (PREIM)}\label{sec:PREIM}
In this section, we first present the main ideas of the PREIM algorithm. Then
we describe one important building block called UPDATE\_EIM. Finally,
using this building block together with the procedure UPDATE\_RB from
Algorithm~\ref{progressive_pod}, we present the PREIM algorithm.

\subsection{Main ideas}
PREIM consists in a progressive construction of the EIM approximation 
and of the RB. 
The key idea is that,
unlike the standard EIM for which HF trajectories are computed for all 
the parameter values in the training set $\mathcal{P}^\mathrm{tr}$ 
(Algorithm~\ref{offline_eim},
line~\ref{off_comp}), PREIM works with an additional training
subset $\mathcal{P}_m^\textsc{HF} \subset \mathcal{P}^\mathrm{tr}$
that is enriched progressively with the iteration index $m$ of PREIM.
The role of $\mathcal{P}_m^\textsc{HF}$ is to collect 
the parameter values for which a HF trajectory has already been computed. 
PREIM is designed so that $\mathrm{Card}(\mathcal{P}_m^\textsc{HF}) \leq m$ for all $m\in\{1,\ldots,M\}$. 
This means that when the final rank-$M$ EIM approximation
has been computed, at most $M$ HF trajectories have been evaluated, 
whence the computational gain with respect to the standard offline stage
provided $M\ll P$.

At the iteration $m\ge1$ of PREIM, the trajectories for all
$\mu \in \mathcal{P}_m^\textsc{HF}$ are HF trajectories, whereas
they are approximated by RB trajectories for all $\mu \in \mathcal{P}^\mathrm{tr}
\setminus \mathcal{P}_m^\textsc{HF}$. The RB functions can be modified
at each iteration $m$ of PREIM; this happens whenever a new value of the parameter is selected 
in the greedy stage of the EIM.
To reflect this, we add a superscript $m$
to the RB trajectories which are now denoted 
$(\hat{u}^{m,k}_\mu)_{k\in \overline{\mathbb K}^{\mathrm{tr}}}$
for all $\mu \in\mathcal{P}^\mathrm{tr}
\setminus \mathcal{P}_m^\textsc{HF}$. It is convenient to
introduce the notation
\begin{equation}\label{preim_func}
\bar u^{m,k}_\mu := \left\{
 \begin{alignedat}{2}
  &u^k_\mu &&\quad\text{if }\mu \in 
 \mathcal P_{m}^\textsc{HF}, \\
  &\hat{u}^{m,k}_\mu &&\quad \text{otherwise},
 \end{alignedat}
 \right .
\end{equation}
and the nonlinear function
\begin{equation} \label{eq:bar_gamma_m}
 \bar \gamma^m(\mu,k,x) := \Gamma(\mu,\bar u_\mu^{m,k}(x)).
\end{equation}
The goal of every PREIM iteration is twofold: 
\begin{enumerate}[(i)]
 \item produce a set of RB functions
$(\theta_n^m)_{1\leq n\leq N^m}$ (the RB functions
and their number depend on $m$); 
\item produce a rank-$m$ approximation of the nonlinear function
$\bar \gamma^m$ defined by~\eqref{eq:bar_gamma_m} in the form
\begin{equation} \label{eq:gamma_m_m}
\bl{
\bar\gamma_{[\mathcal{P}^\textsc{HF}_m,\mathcal{X}_m,\mathcal{Q}_m]}^m(\mu,k,x):=
\sum_{j=1}^m(\bar\varphi^m)_{\mu,j}^{k}\bar q_j(x).
}
\end{equation}
\end{enumerate}
\bl{
The notation $\bar\gamma_{[\mathcal{P}_m^\textsc{HF}, \mathcal{X}_m,\mathcal{Q}_m]}^m$ in~\eqref{eq:gamma_m_m} 
 indicates the data $[\mathcal{P}^\textsc{HF}_m, \mathcal{X}_m,\mathcal{Q}_m]$ that is used to build the 
 approximation of the nonlinearity.} More precisely, this construction uses
the PREIM training set  
$\mathcal{P}_m^\textsc{HF}$,
the sequence of interpolation points 
$\mathcal{X}_m:=(\bar x_i)_{1\le i\le m}$ in $\Omega^{\mathrm{tr}}$ 
(with $\bar x_m$ computed at iteration $m$),
and the sequence of functions $\mathcal{Q}_m:=(\bar q_j)_{1\le j\le m}$ 
defined on $\Omega$  
(with $\bar q_m$ computed at iteration $m$). The progressive construction
of these three ingredients is described below. Then, 
considering the (invertible) lower-triangular 
matrix $\mathbf{\bar B}\in\mathbb R^{m\times m}$
whose last row is calculated using 
$\mathbf {\bar B}_{mj} = \bar q_j(\bar x_m)$ 
for all $j\in\{1,\ldots,m\}$, we compute the real numbers
$(\bar\varphi^m)_{\mu,j}^{k}$ in~\eqref{eq:gamma_m_m} from the relations
\begin{equation}\label{interp}
\sum_{j=1}^m\mathbf{\bar B}_{ij}(\bar\varphi^m)_{\mu,j}^{k}
=\bar\gamma^m(\mu,k,\bar x_i), \qquad \forall  i\in\{1,\ldots,m\},
\end{equation}
for all $(\mu,k) \in\mathcal P\times \overline{\mathbb K}^\mathrm{tr}$.
All the real numbers $(\bar\varphi^m)_{\mu,j}^{k}$ depend on $m$
since the right-hand side of~\eqref{interp} depends on $m$. 

\subsection{The procedure UPDATE\_EIM}

\begin{algorithm}[htb]
\caption{\bl{UPDATE\_EIM\label{improve_eim}}}
\begin{algorithmic}[1]
\vspace{0.2cm}
\Statex {\scshape{\bfseries \underline {Input :}}} $(\theta_n)_{1\leq n \leq N^{m-1}}$, 
$\mathcal{P}_{\mathrm{in}}^\textsc{HF}$, $\mathcal{X}_{m-1}$, 
$\mathcal{Q}_{m-1}$, and $\epsilon_\textsc{eim}$
\State Compute $(\bar{u}^{k}_{\mu})_{(\mu,k)\in \mathcal P^\mathrm{tr}
\times \overline{\mathbb K}^\mathrm{tr}} $ using $(\theta_n)_{1\leq n \leq N^{m-1}}$ \label{line_pr_func}
\State Search $(\mu_m, k_m) 
\in \underset{(\mu',k')\in \mathcal{P}^\mathrm{tr} \times \overline{\mathbb K}^\mathrm{tr}}
{\operatorname{argmax}} \|\Gamma\big(\mu',\bar u^{k'}_{\mu'}(\cdot)\big)
-\bar\gamma_{[\mathcal{P}^\textsc{HF}_{\mathrm{in}},\mathcal{X}_{m-1},\mathcal{Q}_{m-1}]}^{m-1}(\mu',k',\cdot)\|_{\ell^\infty(\Omega^\mathrm{tr})} $ \label{argmax_pr}
 \hfill {\small\verb|based on RB/HF|}
 \State Define $\widetilde{r}_m(\cdot) = \Gamma\big(\mu_m,\bar u^{k_m}_{\mu_m}(\cdot)\big)
-\bar\gamma_{[\mathcal{P}^\textsc{HF}_{\mathrm{in}},\mathcal{X}_{m-1},\mathcal{Q}_{m-1}]}^{m-1}(\mu_m,k_m,\cdot)$. \label{l:widetilde}
\If{$\mu_m \notin \mathcal{P}^\textsc{HF}_{\mathrm{in}}$}\label{9}
\State Compute 
$\mathcal S_{\mathrm{out}} = (u^k_{\mu_m})_{k\in \overline{\mathbb K}^\mathrm{tr}}$ 
and set $\mathcal P^\textsc{HF}_{\mathrm{out}}=\mathcal P^\textsc{HF}_{\mathrm{in}} 
\cup \{ \mu_m\}$
\hfill\hfill {\small\verb|one HF trajectory|} \label{one_HF_trajectory}
\State Search $(\bar\mu_m,\bar k_m)
\in \underset{(\mu',k')\in
{\mathcal{P}_{\mathrm{out}}^\textsc{HF}} \times \overline{\mathbb K}^\mathrm{tr}
}{\operatorname{argmax}}\|\Gamma\big(\mu', u^{k'}_{\mu'}(\cdot)\big)
-\bar\gamma_{[\mathcal{P}^\textsc{HF}_{\mathrm{in}},\mathcal{X}_{m-1},\mathcal{Q}_{m-1}]}^{m-1}(\mu',k',\cdot)\|_{\ell^\infty(\Omega^\mathrm{tr})}$\label{line_upd}
\Else
\State Set $\mathcal S_{\mathrm{out}} = \emptyset$, $\mathcal{P}^\textsc{HF}_{\mathrm{out}} = \mathcal{P}^\textsc{HF}_{\mathrm{in}}$, and $(\bar\mu_m,\bar k_m)= (\mu_m, k_m)$ \label{no_new_trajectory}
\EndIf
\State Define $\bar r_m(\cdot):= \Gamma\big(\bar \mu_m, u^{\bar k_m}
_{\bar \mu_m}(\cdot)\big) 
-\bar\gamma_{[\mathcal{P}^\textsc{HF}_{\mathrm{in}},\mathcal{X}_{m-1},\mathcal{Q}_{m-1}]}^{m-1}(\bar\mu_m,{\bar k}_m,\cdot)$ \label{l:bar}
\If{
$\|\bar r_m\|_{\ell^\infty(\Omega^\mathrm{tr})}<\epsilon_\textsc{eim}$} \label{check_bar_rm}
\State Set ${\rm incr\_rk}$ = FALSE
\State Define $r_m(\cdot)
= \widetilde{r}_m(\cdot)$
\hfill\hfill  {\small\verb|discard the EIM selection|}
\State Set $\mathcal X_{\mathrm{out}} = \mathcal X_{m-1}$ and $\mathcal Q_{\mathrm{out}} = \mathcal Q_{m-1}$
\Else
\State Set ${\rm incr\_rk}$ = TRUE
\State Define $r_m(\cdot)
= {\bar r}_m(\cdot)$
\State Set $\mathcal X_{\mathrm{out}} = (\mathcal X_{m-1},\bar x_m)$ and $\mathcal Q_{\mathrm{out}} = (\mathcal Q_{m-1},\bar q_m)$ with $\bar x_m$, $\bar q_m$ 
as in Algorithm~\ref{offline_eim}
(lines \ref{8} and \ref{rm_xm_eim}).
\EndIf
\State Define $\delta^\textsc{eim}_m =  \|r_m\|_{\ell^\infty(\Omega^\mathrm{tr})}$
\Statex {\scshape{\bfseries \underline {Output :}}}  
${\rm incr\_rk}$, $\mathcal{P}^\textsc{HF}_{\mathrm{out}}$, $\mathcal S_{\mathrm{out}}$, $\mathcal X_{\mathrm{out}}$, 
$\mathcal Q_{\mathrm{out}}$, and $\delta^\textsc{eim}_m$
\vspace{0.2cm}
\end{algorithmic}
\end{algorithm}

An essential building block of PREIM is the procedure UPDATE\_EIM described in Algorithm~\ref{improve_eim}. The input is the RB functions $(\theta_n)_{1\leq n \leq N^{m-1}}$, the triple 
$[\mathcal{P}_{\mathrm{in}}^\textsc{HF},\mathcal{X}_{m-1}, 
\mathcal{Q}_{m-1}]$ describing the current approximation of the nonlinearity (the choice for the indices will be made clearer in the next section, and is not important at this stage), and the threshold $\epsilon_\textsc{eim}$.
The output is the flag ${\rm incr\_rk}$ which indicates whether or not the rank of the EIM approximation has been increased, and if ${\rm incr\_rk}=\mathrm{TRUE}$, the additional output is 
the triple $[\mathcal{P}^\textsc{HF}_{\mathrm{out}},\mathcal X_{\mathrm{out}},\mathcal Q_{\mathrm{out}}]$ to devise the new EIM approximation, possibly a new HF trajectory $\mathcal S_{\mathrm{out}}$, and a measure $\delta^\textsc{eim}_m$ on the EIM error.

First (see line~\ref{argmax_pr}), one selects a new pair $(\mu_m,k_m)\in\mathcal{P^\mathrm{tr}} 
\times \overline{\mathbb K}^\mathrm{tr} $ in a greedy fashion as follows:
 \begin{equation}\label{greedy}
(\mu_m, k_m) 
\in \underset{(\mu',k')\in \mathcal{P}^\mathrm{tr} \times \overline{\mathbb K}^\mathrm{tr}}
{\operatorname{argmax}} \|\Gamma\big(\mu',\bar u^{k'}_{\mu'}(\cdot)\big)
-\bar\gamma_{[\mathcal{P}^\textsc{HF}_{\mathrm{in}},\mathcal{X}_{m-1},\mathcal{Q}_{m-1}]}^{m-1}(\mu',k',\cdot)\|_{\ell^\infty(\Omega^\mathrm{tr})}.
\end{equation}
In~\eqref{greedy}, $\bar u^{k'}_{\mu'}$ is defined as in~\eqref{preim_func} 
using the set $\mathcal{P}^\textsc{HF}_{\mathrm{in}}$. Therefore, the selection criterion~\eqref{greedy} exploits the 
knowledge of the HF trajectory for all the parameter values
in $\mathcal{P}^\textsc{HF}_{\mathrm{in}}$, and otherwise
uses a RB trajectory. This is an important
difference with respect to the standard offline stage.
There are now two possibilities: \textup{(i)} either $\mu_m$ is already in
$\mathcal P_{\mathrm{in}}^\textsc{HF}$; then, no new HF trajectory is 
computed and we set 
$\mathcal P_{\mathrm{out}}^\textsc{HF} := \mathcal P_{\mathrm{in}}^\textsc{HF}$
(line~\ref{no_new_trajectory});
\textup{(ii)} or $\mu_m$ is not in $\mathcal P_{\mathrm{in}}^\textsc{HF}$; 
then we compute
a new HF trajectory for the parameter $\mu_m$ and we set 
$\mathcal P_{\mathrm{out}}^\textsc{HF} := \mathcal P_{\mathrm{in}}^\textsc{HF} 
\cup \{ \mu_m\}$ (line~\ref{one_HF_trajectory}). 
Our numerical experiments reported in
Section~\ref{num_res} below will show that at many iterations
of PREIM, the pair $(\mu_m, k_m)$ selected in~\eqref{greedy}
differs from the previously selected pair by the time index and not 
by the parameter value; this means that for many PREIM iterations,
no additional HF computation is performed.
Nonetheless, in case of non-uniqueness of the maximizer in~\eqref{greedy},
one selects, if possible, a trajectory for which the parameter 
is not already in the set $\mathcal{P}^\textsc{HF}_{\mathrm{in}}$ so as to trigger a computation of a new HF trajectory. 

An additional feature of PREIM is that, whenever
a new HF trajectory is actually computed, one can 
either confirm or update the selected pair $(\mu_m, k_m)$
using the following HF-based re-selection criterion (see line~\ref{line_upd}):
\begin{equation}\label{update}
(\bar\mu_m,\bar k_m)
\in \underset{(\mu',k')\in
{\mathcal{P}_{\mathrm{out}}^\textsc{HF}} \times \overline{\mathbb K}^\mathrm{tr}
}{\operatorname{argmax}}\|\Gamma\big(\mu', u^{k'}_{\mu'}(\cdot)\big)
-\bar\gamma_{[\mathcal{P}^\textsc{HF}_{\mathrm{in}},\mathcal{X}_{m-1},\mathcal{Q}_{m-1}]}^{m-1}(\mu',k',\cdot)\|_{\ell^\infty(\Omega^\mathrm{tr})}.
\end{equation}
We notice that this re-selection criterion only handles HF trajectories since
the parameter values are in $\mathcal P_{\mathrm{out}}^\textsc{HF}$. Moreover, \eqref{update} only requires to probe the values
for $\mu_m$, since the values for the other parameters, which are in 
$\mathcal P_{\mathrm{in}}^\textsc{HF}$, have already been evaluated in~\eqref{greedy}. \bl{Finally, to prevent division by small quantities, the value of the residual $\|\bar r_m\|_{\ell^\infty(\Omega^{\mathrm{tr}})}$ is checked in line~\ref{check_bar_rm}. If this value is too small, the re-selected pair 
$(\bar\mu_m, \bar k_m)$ is rejected and the rank of the EIM approximation is not increased.}

\subsection{The PREIM algorithm}

\begin{algorithm}[htb]
\caption{\bl{Progressive RB-EIM (PREIM)\label{preim}}}
\begin{algorithmic}[1]
\vspace{0.2cm}
\Statex {\scshape{\bfseries \underline {Input :}}} 
$\mathcal{P}^\mathrm{tr}$,
$\overline{\mathbb K}^\mathrm{tr}$,
$\Omega^\mathrm{tr}$,
$ \epsilon_\textsc{pod}>0$, $\epsilon_\textsc{eim}>0$, and $\epsilon_\textsc{rb}>0$
\State Set $m=1$
\State Choose $\mathcal{P}_1^\textsc{HF}\subsetneq\mathcal{P}^\mathrm{tr}$ and
compute $\mathcal S_1 = (u^k_{\mu})_{(\mu,k) \in \mathcal{P}^\textsc{HF}_1\times \overline{\mathbb K}^\mathrm{tr}}$\label{PHF1} 
 \hfill\hfill $J\geq1$ {\small\verb|HF trajectories|}
\State Compute $(\theta_n^1)_{1\leq n\leq N^1}= \mathrm{POD}(\mathcal S_1,
\epsilon_{\textsc{pod}})$.
\State Compute $\hat{\mathbf u}^0\in \mathbb R^{N^1}$, 
$(\mathbf f^k)_{k\in{\mathbb K}^\mathrm{tr}}\in (\mathbb R^{N^1})^K$, 
$\mathbf M\in \mathbb R^{N^1\times N^1}$, and $\mathbf A_0\in \mathbb R^{N^1\times N^1}$
\State Compute $(\mathcal{X}_1,\mathcal{Q}_1,\delta_{1}^\textsc{eim}) = \mathrm{INIT\_EIM}(\mathcal{P}_1^\textsc{HF})$ and $\mathbf C^1 \in \mathbb R^{N^1\times N^1}$
\State Compute
$\delta_1^\textsc{rb}=\operatorname{max}_{\mu\in{\mathcal{P}^\mathrm{tr}}}\Delta_{1}(\mu)$\label{line:zero}
\While{($\delta_{m}^\textsc{eim}>\epsilon_\textsc{eim}$ or
$\delta_{m}^\textsc{rb}>\epsilon_\textsc{rb}$)}\label{while_test}
\State Set $m = m+1$ and $\mathcal{P}_{\mathrm{in}}^\textsc{HF}:=\mathcal{P}_{m-1}^\textsc{HF}$
\State (${\rm incr\_rk}$, $\mathcal{P}_{\mathrm{out}}^\textsc{HF}$, $\mathcal S_{\mathrm{out}}$,
$\mathcal{X}_{\mathrm{out}}$, $\mathcal{Q}_{\mathrm{out}}$,
$\delta_m^\textsc{eim}$) = UPDATE\_EIM
($(\theta_n^{m-1})_{1\leq n\leq N^{m-1}}$, $\mathcal{P}_{\mathrm{in}}^\textsc{HF}$, $\mathcal{X}_{m-1}$,
$\mathcal{Q}_{m-1}$, $\epsilon_\textsc{eim}$)\label{update_eim1}
\While{${\rm incr\_rk}$ = FALSE}
\State $\mathcal{P}_{\mathrm{in}}^\textsc{HF} = \mathcal{P}_{\mathrm{out}}^\textsc{HF}$
\State $(\theta_n^{m-1})_{1\leq n\leq N^{m-1}} = 
\mathrm{UPDATE\_RB}\big((\theta_n^{m-1})_{1\leq n\leq N^{m-1}}, \mathcal S_{\mathrm{out}}, 
\epsilon_\textsc{pod}\big)$ \label{rb_update}
\State (${\rm incr\_rk}$, $\mathcal{P}_{\mathrm{out}}^\textsc{HF}$, $\mathcal S_{\mathrm{out}}$,
$\mathcal{X}_{\mathrm{out}}$, $\mathcal{Q}_{\mathrm{out}}$,
$\delta_m^\textsc{eim}$) = UPDATE\_EIM
($(\theta_n^{m-1})_{1\leq n\leq N^{m-1}}$, $\mathcal{P}_{\mathrm{in}}^\textsc{HF}$, $\mathcal{X}_{m-1}$,
$\mathcal{Q}_{m-1}$, $\epsilon_\textsc{eim}$)\label{update_eim2}
\If {${\rm incr\_rk}$ = TRUE}
\State {Step to line~\ref{rb_fin_update_a}}
\EndIf
\State Compute $\mu_m \in \underset{\mu\in
{\mathcal{P}^\mathrm{tr}}}{\operatorname{argmax}}\ \Delta^{\mathcal X_{\mathrm{out}}, \mathcal Q_{\mathrm{out}}}_{(\theta_n^{m-1})_{1\leq n\leq N^{m-1}}}(\mu)$
\State Compute $\mathcal S_{\mathrm{out}} = (u^k_{\mu_m})_{k\in \overline{\mathbb K}^\mathrm{tr}}$ \label{one_more_HF}
\hfill\hfill {\small\verb|one HF trajectory|}
\EndWhile
\State Set  
$\mathcal{P}_{m}^\textsc{HF} = \mathcal{P}_{\mathrm{out}}^\textsc{HF}$,
$\mathcal S_m = \mathcal S_{\mathrm{out}}$,
$\mathcal X_m = \mathcal X_{\mathrm{out}}$, and
$\mathcal Q_m = \mathcal Q_{\mathrm{out}}$ \label{rb_fin_update_a}
\State Compute $(\theta_n^m)_{1\leq n\leq N^m} = 
\mathrm{UPDATE\_RB}\big((\theta_n^{m-1})_{1\leq n\leq N^{m-1}}, \mathcal S_m, 
\epsilon_\textsc{pod}\big)$  \label{rb_fin_update}
\State Update $\hat{\mathbf u}^0\in \mathbb R^{N^m}$,
$(\mathbf f^k)_{k\in\mathbb K^\mathrm{tr}}\in (\mathbb R^{N^m})^K$, and the matrices
$\mathbf M$, $\mathbf A_0$, $(\mathbf C^j)_{1\leq j \leq m}$ in $\mathbb R^{N^m\times N^m}$\label{mat_upd}
\State Compute $\delta_m^\textsc{rb}=
\underset{\mu\in{\mathcal{P}^\mathrm{tr}}}{\operatorname{max}}\Delta^{\mathcal X_{m}, \mathcal Q_{m}}_{(\theta_n^{m})_{1\leq n\leq N^{m}}}(\mu)$
\EndWhile
\State Set $M:=m$
\Statex {\scshape{\bfseries \underline {Output :}}}  $(\theta_n)_{1\leq n\leq {N^M}}$, $\hat{\mathbf u}^0$,
$(\mathbf f^k)_{k\in{\mathbb K}^\mathrm{tr}}$, $\mathbf M$, $\mathbf A_0$, 
$\mathcal X_M$, $\mathcal Q_M$, and $(\mathbf C^j)_{1\leq j \leq M}$
\vspace{0.2cm}
\end{algorithmic}
\end{algorithm}

We are now ready to describe the PREIM procedure, see Algorithm~\ref{preim}.
PREIM is an iterative method that builds progressively the RB and the EIM approximation. The iteration is controlled by three tolerances: 
$ \epsilon_\textsc{pod}>0$ which is used in the progressive increment of the RB, $\epsilon_\textsc{eim}>0$ which is used to check the quality of the EIM approximation, and $\epsilon_\textsc{rb}>0$ which is used to check the quality of the RB. The termination criterion involves the quality of both the EIM and the RB approximations, see line~\ref{while_test}. Note that this is the same criterion as in the standard RB-EIM approach.

\bl{Within each PREIM iteration, the two previously-described procedures UPDATE\_EIM and UPDATE\_RB are called. First, one attempts to improve the EIM approximation (line~\ref{update_eim1}). 
If this is successful (i.e., if $\text{incr\_rk}=\text{TRUE}$), the RB is updated by using the possibly new HF trajectory $\mathcal S_m$ (line~\ref{rb_fin_update}). 
Otherwise (i.e., if $\text{incr\_rk}=\text{FALSE}$), the RB is possibly updated (line~\ref{rb_update}) and a new improvement of the EIM is attempted (line~\ref{update_eim2}). 
In general, the RB is improved because a new HF trajectory has been computed. Whenever this is not the case, a new HF trajectory is anyway computed in line~\ref{one_more_HF} so as to steer the progress in the iterations.} The 
choice of this new HF trajectory can be driven by a standard greedy algorithm based on the use of a classical a posteriori error estimator. More precisely, for a given reduced basis $(\theta_n)_{1\leq n \leq N}$ and 
given sets of training points $\mathcal X$ and functions $\mathcal{Q}$ used for the current EIM approximation of the nonlinearity, the associated a posteriori error estimator for a given value of the parameter $\mu \in \mathcal P$ is denoted by 
$
\Delta^{\mathcal X, \mathcal Q}_{(\theta_n)_{1\leq n\leq N}}(\mu)
$.
Finally, we observe that
the reduced matrices and vectors in line~\ref{mat_upd} of Algorithm~\ref{preim} need to be updated
since these quantities depend on the RB functions which can change at 
every iteration.

Let us now discuss the initialization of PREIM. 
\bl{In line~\ref{PHF1}, one can choose an initial PREIM training set $\mathcal P^{\mathrm{HF}}_1$ composed of a single parameter, as is often the case with greedy algorithms. Although the nonlinearity may not be well-described initially, one can expect that the description will improve progressively. Still, to allow for more robustness in the initialization, one can consider an initial set $\mathcal P^{\mathrm{HF}}_1$ composed of several parameters.}
One can then compute the HF trajectories for all $\mu\in \mathcal P^{\mathrm{HF}}_1$ and compress them using the POD procedure with threshold $\epsilon_\textsc{pod}$ (if $\mathcal P^{\mathrm{HF}}_1$ contains more than one value, a progressive version is used). Finally, one selects
\begin{equation} 
(\mu_1, k_1) 
\in \underset{(\mu',k')\in \mathcal{P}_1^\textsc{HF} \times \overline{\mathbb K}^\mathrm{tr}}
{\operatorname{argmax}} \|\Gamma\big(\mu', u^{k'}_{\mu'}(\cdot)\big)\|
_{\ell^\infty(\Omega^\mathrm{tr})},
\end{equation}
one defines $r_1(\cdot)= \Gamma(\mu_1,u^{k_1}_{\mu_1}(\cdot))$
and computes $\mathcal X_1=(\bar x_1)$, $\mathcal Q_1=(\bar q_1)$ (as in the standard EIM procedure), 
and one sets $\delta_1^\textsc{eim} = \|r_1\|_{\ell^\infty(\Omega^\mathrm{tr})}$.

\begin{remark}[PREIM-NR and U-SER variants] \label{rem:NR_USER}
We can consider two variants in the procedure UPDATE\_EIM
(Algorithm~\ref{improve_eim}) and therefore in PREIM.
A first variant consists in skipping the re-selection 
step in line~\ref{line_upd} of Algorithm~\ref{improve_eim}. 
This variant, which we call PREIM-NR (for `no re-selection'),
will be tested numerically in the next section so as to highlight
the actual benefits brought by the re-selection.
A second variant is to 
replace $\bar u^{m,k}_\mu$ with 
$\hat u^{m,k}_\mu$ in lines~\ref{line_pr_func}
and~\ref{argmax_pr} of Algorithm~\ref{improve_eim}, and to skip the re-selection
step in line~\ref{line_upd}. 
We call this variant U-SER since it can be viewed as an extension of SER~\cite{ser15} to the unsteady setting. 
The crucial difference 
between PREIM-NR and U-SER is that U-SER uses RB trajectories to 
compute the space-dependent functions in the EIM approximation whereas
PREIM-NR uses HF trajectories.
\end{remark}

 \section{Numerical results}\label{num_res}
\begin{figure}[htb]
\includegraphics[scale=0.3]{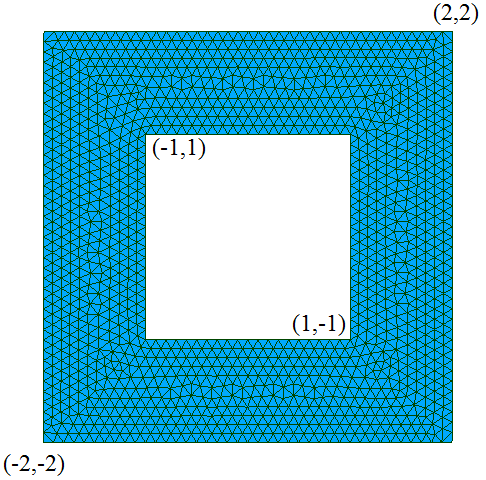} 
\centering 
\caption{Computational domain: perforated plate.}
\label{plate}
\end{figure}

In this section, we illustrate the above developments by numerical 
examples related to transient heat transfer problems with two different types
of nonlinearities. The first example uses a nonlinearity on the solution
whereas the second investigates a nonlinearity on its partial derivatives.
Our goal is to illustrate the computational performance of PREIM
and to compare it to the standard EIM approach described in Section~\ref{off}
and to the variants PREIM-NR and U-SER described in Remark~\ref{rem:NR_USER}. 
We consider a two-dimensional setting based on 
the perforated plate illustrated in Figure~\ref{plate} with
$\Omega=(-2,2)^2\backslash[-1,1]^2\subset \mathbb R^2$.
HF trajectories are computed using a Finite Element subspace
$X\subset Y= H^1(\Omega)$
consisting of continuous, piecewise affine functions. 
The HF computations use the industrial software
\texttt{code\_aster}~\cite{aster} for the first test case and \texttt{FreeFem++}~\cite{freefem} for the second test case, 
whereas the reduced-order modeling algorithms have 
been developed in Python.
In both test cases, the dominant error component turns out to be
the one resulting
from the approximation of the nonlinearity, rather than the one resulting from
the RB. For this reason, PREIM has been run using only the 
stopping criterion $\delta_{m}^\textsc{eim}>\epsilon_\textsc{eim}$ in line~\ref{while_test} of Algorithm~\ref{preim}.

\subsection{Test case (a): Nonlinearity on the solution}

\begin{figure}[htb]
\includegraphics[scale=0.25]{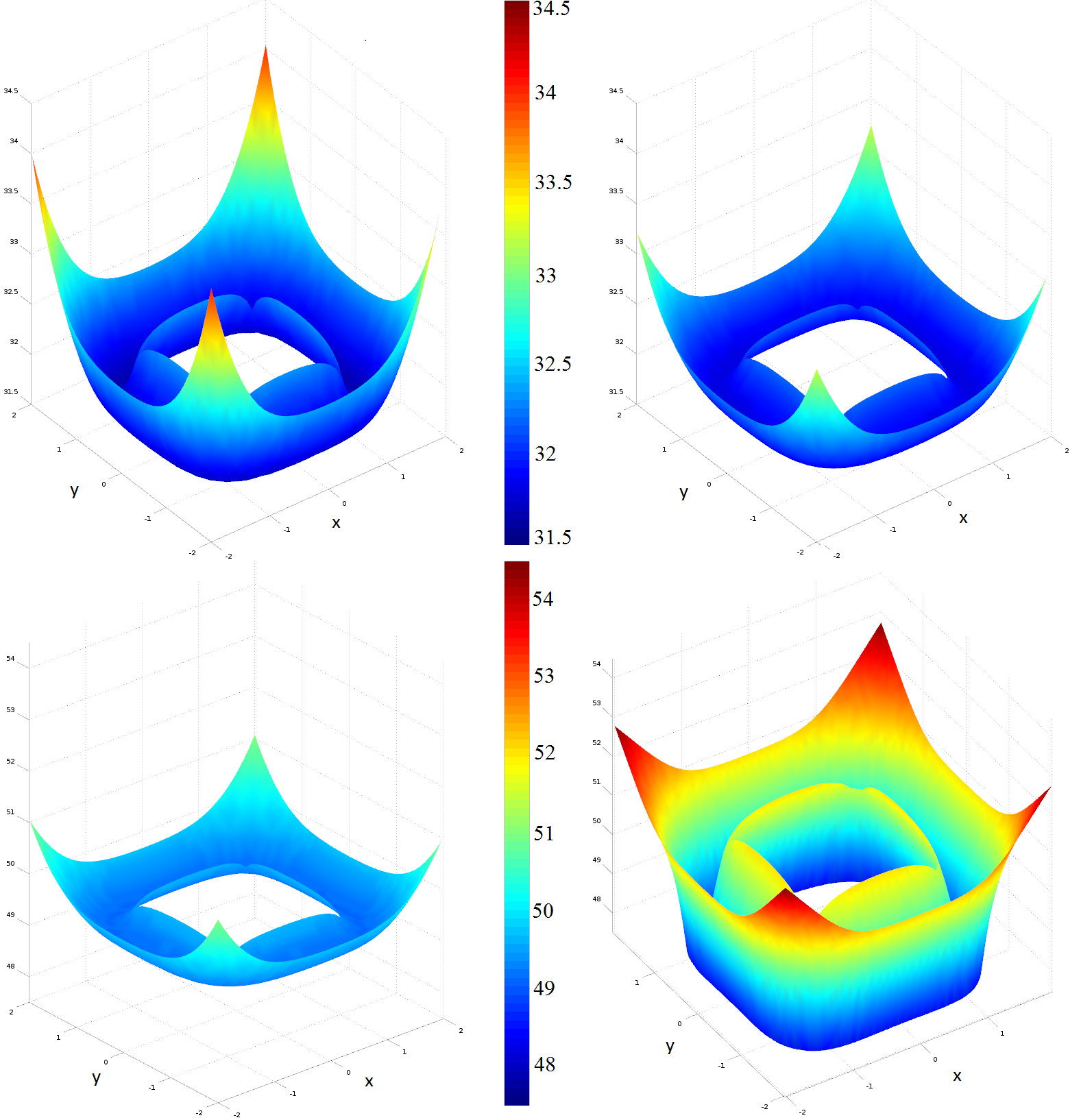} 
\centering 
\caption{Test case (a): HF solutions for the parameter values
$\mu = 1$ (left) and $\mu = 18$ (right) at $t=2$~s (top) 
and $t=5$~s (bottom).}
\label{bulk}
\end{figure}

We consider the nonlinear parabolic problem~\eqref{HTE} with
the nonlinear function $\Gamma(\mu,v) := 
\sin\Big(\frac{2\pi\mu}{20}\big(\frac{v-u_0}{u_{m}-u_0}\big)^2\Big)$, 
with $u_0=293$~K (20~$^o$C) and $u_m=323$~K (50~$^o$C). 
We define $\kappa_0 = 1.05$ m$^{2}\cdot$K$^{-2}\cdot$s$^{-1}$ 
and $\phi_e = 3$ K$\cdot$m$\cdot$s$^{-1}$ 
(these units result from our normalization by the density times the heat capacity).
For space discretization, we use a mesh containing $\mathcal N = 1438$ nodes (see Figure~\ref{plate}).
Regarding time discretization, we consider the time interval $I = [0,5]$,
the set of discrete times nodes $\mathbb K^\mathrm{tr} = \{1,\cdots,50\}$,
and a constant time step $\Delta t^k=0.1$~s 
for all $k \in \mathbb K^\mathrm{tr}$.
Finally, we consider the parameter interval $\mathcal P = [1,20]$, 
the training set $\mathcal P^\mathrm{tr} = \{1,\cdots,20\}$, and 
\bl{we use the larger set $\{0.25i\ |\ 0\leq i\leq 80\}$ 
to verify our numerical results.}
In Figure \ref{bulk}, we show the HF temperature profiles
over the perforated plate at
two different times and for two different parameter values.
We can see that, as the simulation time increases, 
the temperature is, overall,
higher for larger values of the parameter $\mu$ 
than for smaller values.
Also, for larger values of 
$\mu$, the temperature variation tends to be less uniform over the plate 
than for smaller values of $\mu$. 
\begin{table}[htb]
	\centering
		\begin{tabular}{|c||c|c|c|c|c|c|c|}
		\hline
      $m$ &$1$ & $2$ & $6$ &$14$ & $15$ &$20$ &$25$  \\ \hline
      $\|r_m\|_{\ell^\infty(\Omega^\mathrm{tr})}$ &
      $2.0$& $8.1\mathrm{E}{-1}$ & $1.1\mathrm{E}{-1}$ & $5.2\mathrm{E}{-3}$
      & $2.6\mathrm{E}{-3}$ & $1.1\mathrm{E}{-3}$& $1.6\mathrm{E}{-4}$ \\ \hline
      	\end{tabular}
	\caption{Test case (a): Evolution of the standard EIM error. 
	$m$ is the rank of the EIM approximation.}
	\label{eim_tb}
\end{table}

During the standard offline stage, we perform $P = 20$ 
HF computations. Knowing that $K = 50$, 
the set $\mathcal S$~(Algorithm \ref{offline_eim}, line \ref{off_comp}) 
contains $1020$ fields, each consisting of $\mathcal N=1438$ nodal values. 
Applying the POD in a progressive manner (see Algorithm~\ref{offline_rb} with the parameters enumerated using increasing values) based on the $H^1$-norm
and a truncation threshold $\epsilon_\textsc{pod}=10^{-3}$,
we obtain $N=6$ RB functions. 
Afterwards, we perform the standard EIM algorithm
whose convergence is reported in Table~\ref{eim_tb}. 
For $\epsilon_\textsc{eim}=5\cdot10^{-2}$, the final rank of the EIM approximation is $M=8$, whereas for $\epsilon_\textsc{eim}=10^{-3}$, the final rank of the EIM approximation is $M=15$.

\begin{table}[htb]
	\centering
		\begin{tabular}{|c|c||c|c|c|c|c|c|c|c|c|c|c|c|c|}
		\hline
      \multicolumn{2}{|l||}{$\qquad m$} & $1$&$2$&$3$&$4$&$5$&$6$&$7$&$8$&$9$&$10$
                &$11$&$12$&$13$\\
		\hline\hline
      \multirow{3}*{PREIM}&$\bar\mu$  &$1$ & $20$&
                $20$& $20$& $20$ &$20$& $20$& $20$&          
                $16$ & $20$& $20$& $18$& $20$\\
	        \cline{2-15}
	        &$\mu$  & \cellcolor{gray}$1$ & \cellcolor{gray}$20$&
                $20$& $20$& $20$ &\cellcolor{gray}$18$&$20$&  $20$&          
                \cellcolor{gray}$16$ & $20$& $20$& $18$& $20$\\
	        \cline{2-15}
		&$k$ &$50$ & $45$& $48$& $50$& 
		$43$& $42$& $39$& $46$
		& $50$& $49$& $33$& $50$& $47$\\
		\hline
		\end{tabular}
	\caption{Test case (a):
	Selected parameters and time nodes in PREIM. 
	The gray cells correspond to a new parameter selection and, therefore, 
	to a new HF computation.}
	\label{mu_tb}
\end{table}

We now investigate PREIM, which we first run with thresholds
$\epsilon_\textsc{pod}=10^{-3}$ and $\epsilon_\textsc{eim}=5\cdot10^{-2}$.
Table~\ref{mu_tb} shows the selected parameters and discrete time nodes
at each stage of PREIM.
We can make two important observations from this table.
First, after 13 iterations, PREIM has only selected four different parameter
values, and has therefore computed only four HF trajectories (the iterations for which a new parameter value is selected are indicated in gray in Table~\ref{mu_tb}). In the other $9$ out of the 13 iterations, a different time snapshot of an already existing 
HF trajectory has been selected. Second, by comparing the lines in Table~\ref{mu_tb} related to $\mu$ and $\bar\mu$, we can see that a parameter re-selection happened at iteration $m=7$. 

\begin{figure}[htb]
\includegraphics[scale=0.35]{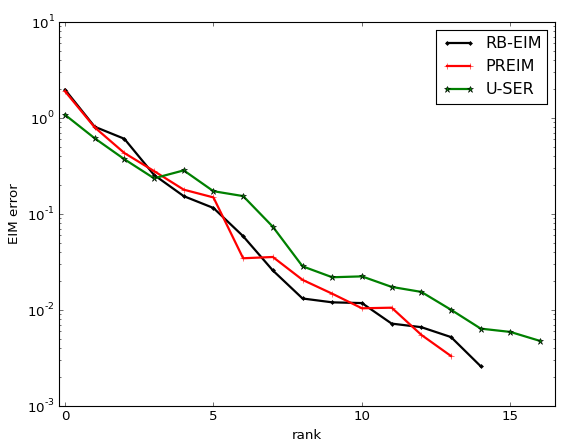}
\includegraphics[scale=0.47]{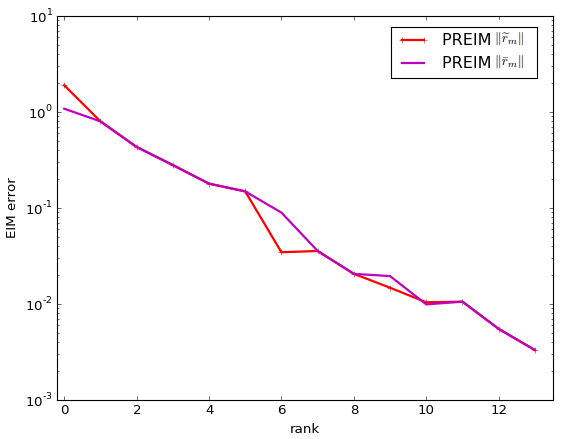}
\centering 
\caption{Test case (a): EIM approximation error
as a function of $m$
for $\epsilon_\textsc{pod}=10^{-3}$ and $\epsilon_\textsc{eim}=5\cdot10^{-2}$.
Left: Errors for the standard RB-EIM procedure, PREIM, and U-SER.
Right: Errors $\|\widetilde r_m\|_{\ell^\infty(\Omega^{\mathrm{tr}})}$ and $\|\bar r_m\|_{\ell^\infty(\Omega^{\mathrm{tr}})}$ for PREIM.
}
\label{sg_vals}
\end{figure}

The left panel of Figure~\ref{sg_vals} displays the error on 
the approximation of the nonlinear function $\Gamma$ for the 
standard RB-EIM procedure and for  
PREIM as a function of the iteration number $m$ (the additional curve concerning U-SER will be commented afterwards), i.e., we plot $\|\bar r_ {m}\|_{\ell^\infty(\Omega^\mathrm{tr})}$
(line~\ref{l:widetilde} of Algorithm~\ref{improve_eim}) and
$\|\widetilde{r}_{m}\|_{\ell^\infty(\Omega^\mathrm{tr})}$ (line~\ref{l:bar} of Algorithm~\ref{improve_eim})
as a function of $m$, see~\eqref{eq:gamma_m_m}. 
We can see that the quality of the approximation of the nonlinearity is almost the same for PREIM as for the standard RB-EIM procedure; yet, the former achieves this accuracy by computing 20\% of the HF trajectories computed by the latter (4 instead of 20 HF trajectories). The right panel of Figure~\ref{sg_vals} shows
the values of $\|\widetilde{r}_m\|_{\ell^\infty(\Omega^\mathrm{tr})}$ and $\|\bar r_m\|_{\ell^\infty(\Omega^\mathrm{tr})}$ as a function of $m$. The two quantities differ when the parameter $\mu_m$ in 
line~\ref{argmax_pr} of Algorithm~\ref{improve_eim} is not in the set $\mathcal P^{\mathrm{HF}}_{m-1}$ so that $\|\widetilde{r}_m\|_{\ell^\infty(\Omega^\mathrm{tr})}$ is computed using a RB approximation whereas $\|\bar r_m\|_{\ell^\infty(\Omega^\mathrm{tr})}$ results from a HF trajectory. Discarding the initialization, 
this happens for $m\in\{6,9,10\}$. The fact that $\|\widetilde{r}_m\|_{\ell^\infty(\Omega^\mathrm{tr})}$ and $\|\bar r_m\|_{\ell^\infty(\Omega^\mathrm{tr})}$ take rather close values indicates that the RB provides an accurate approximation of the HF trajectory.

\begin{figure}[htb]
\includegraphics[scale=0.35]{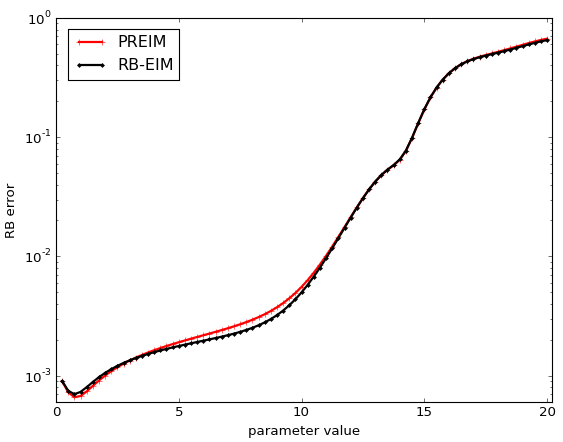}
\includegraphics[scale=0.35]{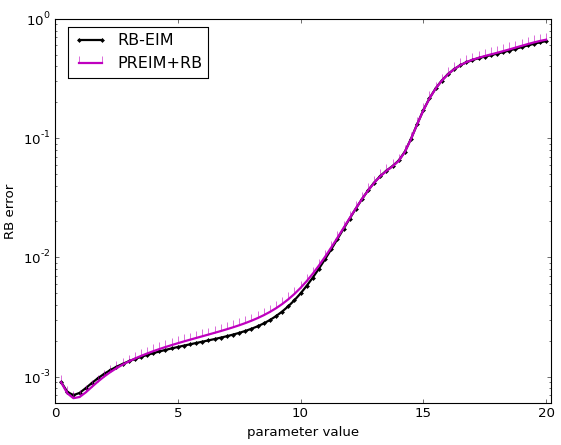}
\centering 
\caption{Test case (a):
RB approximation error $\|u_\mu -\hat u_\mu\|
_{\ell^2(I^\mathrm{tr};H^1(\Omega^\mathrm{tr}))}$ for 
$\epsilon_\textsc{pod}=10^{-3}$ and $\epsilon_\textsc{eim}=5\cdot10^{-2}$.}
 \label{rb_err}
\end{figure}
\begin{figure}
\includegraphics[scale=0.35]{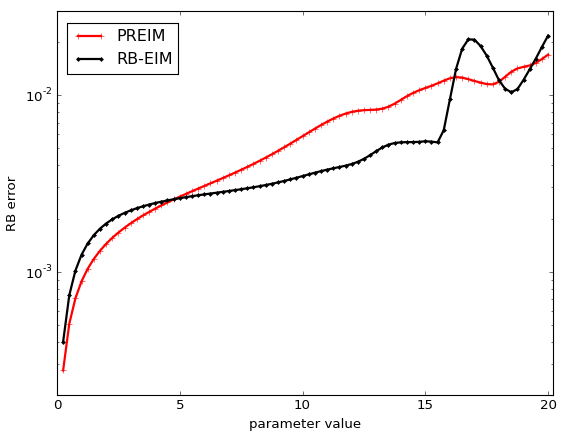}
\includegraphics[scale=0.35]{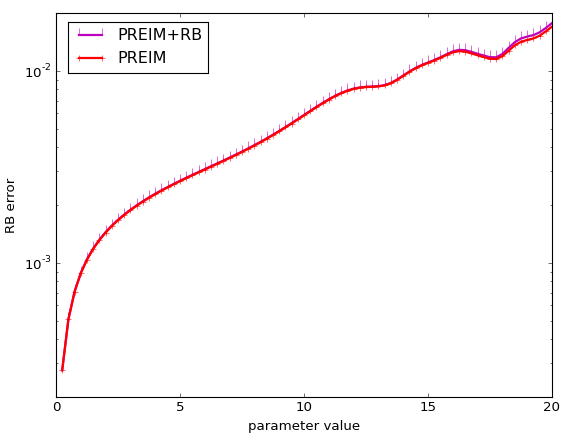}
\centering 
\caption{Test case (a):
RB approximation error $\|u_\mu -\hat u_\mu\|
_{\ell^2(I^\mathrm{tr};H^1(\Omega^\mathrm{tr}))}$ for $\epsilon_\textsc{pod}=10^{-5}$
and $\epsilon_\textsc{eim}=5\cdot10^{-3}$.}
\label{smaller_tol_blk}
\end{figure}
The left panel of Figure~\ref{rb_err} compares the space-time
errors \bl{(measured using the $\ell^2$-norm in time and the
$H^1$-norm in space)} on the trajectories produced by the standard RB-EIM and the PREIM procedures for the whole parameter range. The error is generically denoted $\|u_\mu -\hat u_\mu\|_{\ell^2(I^\mathrm{tr};H^1(\Omega^\mathrm{tr}))}$. We observe an excellent agreement over the whole parameter range. In the right panel of Figure~\ref{rb_err}, we also consider the space-time errors on the trajectories produced using the approximation of the nonlinearity resulting from PREIM with the RB resulting from the standard algorithm. We do not observe any significant change with respect to the left panel, which indicates that the dominant error component is that associated with the approximation of the nonlinearity.
We consider the tighter couple of thresholds $\epsilon_\textsc{pod}=10^{-5}$
and $\epsilon_\textsc{eim}=5\cdot10^{-3}$ in Figure~\ref{smaller_tol_blk}.
Here, we can observe some differences in the errors produced by the standard RB-EIM and PREIM procedures, although both errors remain comparable and reach similar maximum values over the parameter range. While the standard procedure is slightly more accurate for most parameter values, the conclusion is reversed for some other values. Moreover, the curves on the right panel of Figure~\ref{smaller_tol_blk} corroborate the fact that once again, the dominant error component is that associated with the approximation of the nonlinearity.

\begin{figure}
\includegraphics[scale=0.35]{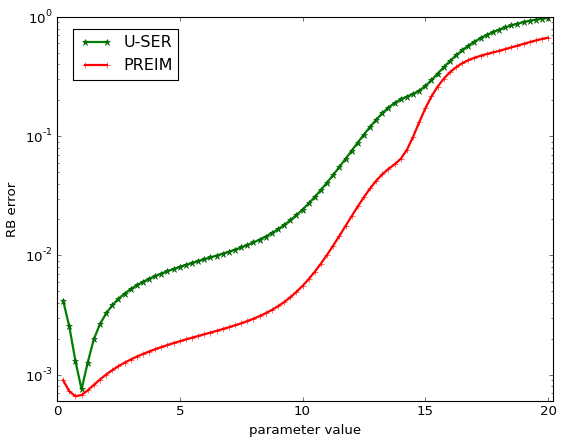}
\centering 
\caption{Test case (a):
RB approximation error $\|u_\mu -\hat u_\mu\|
_{\ell^2(I^\mathrm{tr};H^1(\Omega^\mathrm{tr}))}$ for 
$\epsilon_\textsc{pod}=10^{-3}$ and $\epsilon_\textsc{eim}=5\cdot10^{-2}$.}
 \label{rb+preim}
\end{figure}

\begin{table}[htb]
	\centering
		\begin{tabular}{|c|c||c|c|c|c|c|c|c|c|c|c|c|c|c|}
		\hline
      \multicolumn{2}{|l||}{$\qquad m$} & $1$&$2$&$3$&$4$&$5$&$6$&$7$&$8$&$9$&$10$
                &$11$&$12$&$13$\\
		\hline\hline
	\multirow{2}*{U-SER}&$\mu$ & \cellcolor{gray}$1$ & \cellcolor{gray}$20$& $20$& 
	        $20$& \cellcolor{gray}$16$& $20$& \cellcolor{gray}$19$& $20$
		& $20$& $19$& \cellcolor{gray}$17$& $20$& $19$ \\  
		\cline{2-15}
		&$k$ &$50$ & $49$& $50$& $46$& $42$& 
		$49$& $44$& $39$
		& $50$& $49$& $48$& $47$& $50$\\
		\hline
		\hline
	\multirow{2}*{PREIM-NR}&$\mu$ & \cellcolor{gray}$1$ & \cellcolor{gray}$20$& $20$& 
	        $20$& $20$& \cellcolor{gray}$16$& $20$& $20$
		& $20$& $20$& $20$& \cellcolor{gray}$17$& \cellcolor{gray}$19$ \\  
		\cline{2-15}
		&$k$ &$50$ & $47$& $50$& $46$& $42$& 
		$49$& $48$& $46$
		& $39$& $50$& $45$& $50$& $50$\\
		\hline
		\end{tabular}
	\caption{Test case (a):
	Selected parameters and time nodes in U-SER and PREIM-NR. 
	The gray cells correspond to a new parameter selection and, therefore, 
	to a new HF computation.}
	\label{mu_tb_2}
\end{table}
Let us further explore the PREIM algorithm by comparing it to its variants
U-SER and PREIM-NR introduced in Remark~\ref{rem:NR_USER}.
Table~\ref{mu_tb_2} reports the selected parameters and time nodes in U-SER and PREIM-NR (compare with Table~\ref{mu_tb} for PREIM). 
Both U-SER and PREIM-NR need to
compute five HF trajectories, which is only 25\% of those needed with the
standard RB-EIM procedure, but this is still one more HF trajectory than 
with PREIM.
One difference between U-SER and PREIM-NR is that new 
parameters are selected earlier with U-SER. Interestingly, after
13 iterations, U-SER and PREIM-NR have selected the same five parameters. 
Another interesting observation is that U-SER actually
selects the same couple $(\mu,k)$
twice (this happens for $m=2$ and $m=6$); this can be interpreted by
observing that owing to the improvement of the RB using HF trajectories 
between iterations $m=2$ and $m=6$, the algorithm detects the need to 
improve the approximation of the nonlinearity by using the same pair 
$(\mu,k)$. The same observation can be made for PREIM-NR 
(this happens for $m=4$ and $m=8$). We emphasize that re-selecting the same
pair $(\mu,k)$ is not possible within PREIM since the selection is based on HF trajectories. The left panel of Figure~\ref{sg_vals} displays the error on 
the approximation of the nonlinear function $\Gamma$ obtained with U-SER
and compares it to the error obtained with the standard RB-EIM and PREIM procedures that were already discussed. The U-SER error is evaluated as
$\sup_{(\mu,k)\in \mathcal P^\mathrm{tr}\times \overline{\mathbb K}^\mathrm{tr}}\|\Gamma(\mu,\hat u^k_\mu(\cdot))-\overline\gamma_{[\mathcal{P}_m^\textsc{HF}, \mathcal{X}_m,\mathcal{Q}_m]}^m(\mu,k,\cdot)\|_{\ell^\infty(\Omega^\mathrm{tr})}$. We observe that the approximation of the nonlinearity is somewhat less sharp with U-SER than with PREIM.
Figure~\ref{rb+preim} reports the space-time errors
(measured using the $\ell^2$-norm in time and the
$H^1$-norm in space) on the trajectories produced by PREIM and U-SER for the whole parameter range. We observe that the U-SER error is always larger, sometimes up to a factor of five, but for the larger parameter values which produce the larger errors, the quality of the results produced by PREIM and U-SER remains comparable. 

\begin{table}[htb]
	\centering
		\begin{tabular}{|c||r||r||r||}
		\hline
      $ $ & RB-EIM & PREIM & U-SER \\ \hline
      HF computations & $99\%$ & $20.0\%$ & $25.0\%$ \\ \hline
      greedy runtime & $1\%$ & $1.5\%$ & $2.3\%$ \\ \hline
      Total runtime & $100\%$ & $21.5\%$ & $26.3\%$ \\ \hline
		\end{tabular}
	\caption{Test case (a): Runtime measurements.}
	\label{runtime}
\end{table}

\bl{Finally, we provide an assessment of the runtimes in Table~\ref{runtime}.
We can see that for the standard RB-EIM procedure, the computation of the HF trajectories dominates the cost of the offline phase. For both PREIM and U-SER, the cost of these HF computations is substantially reduced. At the same time, the cost of the greedy algorithm (which includes the construction of the EIM and of the RB) is increased by 50\% with respect to the standard RB-EIM procedure. However, the impact on the total runtime is marginal.
}

\FloatBarrier

\subsection{Test case (b): \bl{Nonlinearity on the partial derivatives}}

\begin{figure}[htb]
\includegraphics[scale=0.27]{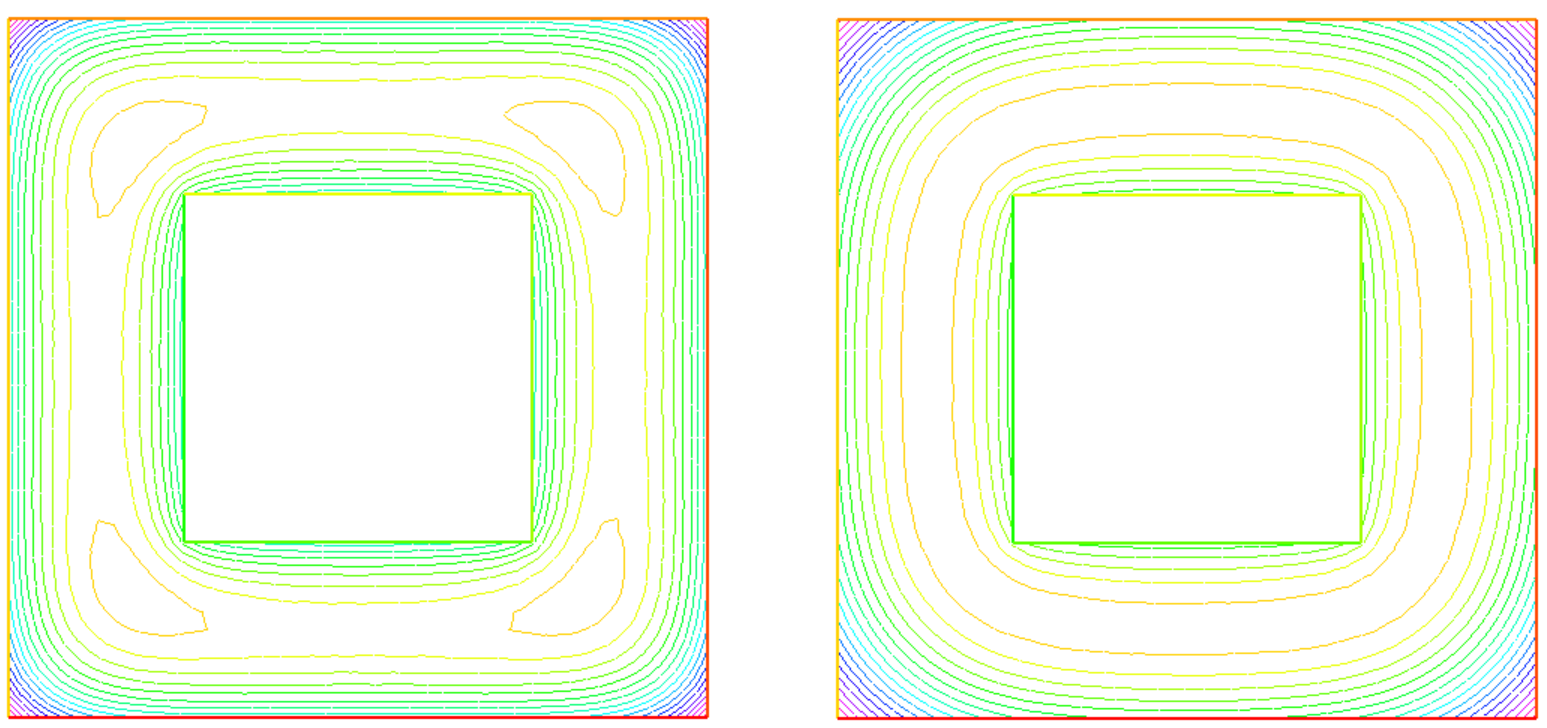}
\centering
\caption{Test case (b): HF solutions for the parameter values $\mu=1$ 
at $t = 1$~s (left, values from 20.2 to 22.1) and at $t = 2.5$~s (right, values from 34.5 to 37.3).}
\label{isoval}
\end{figure}

\bl{We consider the nonlinear parabolic problem~\eqref{HTE} with
the nonlinear function $\Gamma(\mu,u) := 
\sin\Big(\omega\mu\big(({\frac{\partial u}{\partial x}})^2+({\frac{\partial u}{\partial y}})^2\big)\Big)^2$, where $\omega = 6.25\cdot10^{-3}$.
We define $u_0=293$~K (20~$^o$C), $\kappa_0 = 1$ m$^{2}\cdot$K$^{-2}\cdot$s$^{-1}$ and $\phi_e = 3$ K$\cdot$m$\cdot$s$^{-1}$ 
(these units result from our normalization by the density times the heat capacity). 
For the space discretization, we use a mesh containing $\mathcal N = 1429$ nodes.
Regarding time discretization, we consider the time interval $I = [0,2.5]$,
the set of discrete times nodes $\mathbb K^\mathrm{tr} = \{1,\ldots,50\}$,
and a constant time step $\Delta t^k=0.05$~s 
for all $k \in \mathbb K^\mathrm{tr}$.
Finally, we consider the parameter interval $\mathcal P = [1,40]$ 
and the training set $\mathcal P^\mathrm{tr} = \{1,\ldots,40\}$.
In Figure~\ref{isoval}, we show the temperature isovalues
over the perforated plate at two different times for $\mu=1$. 
We can observe different boundary layers 
depending on the time (the same observation can be made by varying the parameter value)}.
\begin{table}[h!]
	\centering
		\begin{tabular}{|c||c|c|c|c|c|c|c|c|c|c|c|c|c|c|}
		\hline
      $ p$ &$1$&$2$&$3$&$8$&$20$&$23$&$24$&$26$&$32$& $33$ & $36$ & $37$& $39$&$40$ \\ \hline
      RB size & $3$ & $4$ & $5$& $6$ & $7$ & $8$& $9$ & $10$ & $11$ & $12$& $13$ & $14$ & $15$&$15$  \\ \hline
		\end{tabular}
	\caption{Test case (b): Size of the reduced basis in the standard algorithm with $ \epsilon_\textsc{pod} = 5\cdot10^{-2}$.}
	\label{b:rb_sz}
\end{table}
\begin{table}[htb]
	\centering
		\begin{tabular}{|c||c|c|c|c|c|c|c|c|c|c|}
		\hline
      $ m$ &$2$& $10$ &$13$&$20$&$30$&$36$&$37$&$79$&$96$&$144$ \\ \hline
      $\|r_m\|_{\ell^\infty(\Omega^{\rm tr})}$ & $1.6$& $1.3$ &  $9.7{\rm E}{-1}$ & $4.7{\rm E}{-1}$& 
      $1.7{\rm E}{-1}$ & $1.2{\rm E}{-1}$
      & $8.0{\rm E}{-2}$ & $9.1{\rm E}{-3}$ & $4.6{\rm E}{-3}$ & $9.4{\rm E}{-4}$ \\ \hline
		\end{tabular}
	\caption{Test case (b): Evolution of the standard EIM error. 
	$m$ is the rank of the EIM approximation 
	and $\|r_m\|_{\ell^\infty(\Omega^\mathrm{tr})}$ 
	is the residual norm in~\eqref{resid}.}
	\label{b:eim}
\end{table}

\bl{During the standard offline stage, we perform $P = 40$ 
HF computations. Knowing that $K = 50$, 
the set $\mathcal S$ (Algorithm \ref{offline_eim}, line \ref{off_comp}) 
contains $2040$ fields, each consisting of $\mathcal N=1429$ nodal values. 
Applying Algorithm~\ref{offline_rb} based on the $H^1$-norm,
a truncation threshold $\epsilon_\textsc{pod}=5\cdot10^{-2}$,
and parameters enumerated with increasing values,
we obtain $N=15$ RB functions. Table~\ref{b:rb_sz} shows the dimension of the RB
space as a function of the enumeration index $p$.
Table~\ref{b:eim} shows the evolution of the error on the nonlinearity 
within the standard EIM.
The fact that the nonlinearity depends on the partial derivatives of the solution challenges the EIM; indeed, the error decay is not as fast as 
in the previous test case. This observation is corroborated by the fact that the functions $(q_j)_{1\le j \le M}$ all look quite different.}
%
%

\begin{table}[htb]
\centering\begin{tabular}{l}
		\begin{tabular}{|c||c|c|c|c|c|c|c|c|c|}
		\hline
                $m$ & $1$&$2$&$3$&$4$&$5$&$6$&$7$&$8$&$9$\\
		\hline\hline
                $\bar\mu$ & $21$ &$8$ & $21$&
                $8$& $21$& $21$ &$21$& $8$&          
                $9$ \\
	        \cline{1-10}
	        $\mu$&\cellcolor{gray}$21$ &\cellcolor{gray}$8$ & $21$&
                $8$& $21$& $21$ &$21$& $8$&\cellcolor{gray}$9$ \\
	        \cline{1-10}
		$k$ &$2$ & $5$& $3$& $2$& 
		$50$& $4$& $49$& $3$
		&$4$\\
		\hline
		\end{tabular}
\\[10mm]
		\begin{tabular}{|c||c|c|c|c|c|c|c|c|c|c|c|c|c|c|c|c|c|c|c|c|}
		\hline
                $m$ & $1$&$2$&$3$&$4$&$5$&$6$&$7$&$8$&$9$&$10$
                &$11$&$12$&$13$&$14$&$15$&$16$&$17$&$18$&$19$&$20$\\
		\hline\hline
                $\bar\mu$  &$21$ & $21$&
                $21$& $8$& $21$ &$21$& $21$& $8$&          
                $9$ & $21$& $9$& $21$& $9$& $9$ & $9$
 		& $6$& $21$& $21$ & $40$&$40$ \\
	        \cline{1-21}
                $\mu$&  \cellcolor{gray}$21$ &\cellcolor{gray}$8$ & $21$&
                $8$& $21$& $21$ &$21$& $8$&          
                \cellcolor{gray}$9$ & $21$& $9$& \cellcolor{gray}$7$&
                \cellcolor{gray}$6$& $9$ & $9$
 		& \cellcolor{gray}$5$& \cellcolor{gray}$4$& \cellcolor{gray}$3$ & \cellcolor{gray}$40$&$40$ \\
	        \cline{1-21}
		$k$ &$2$ & $5$& $3$& $2$& 
		$50$& $4$& $49$& $3$
		&$4$& $10$& $50$& $25$& $49$
		& $5$& $10$& $4$& $6$& $9$ & $15$& $40$\\
		\hline
		\end{tabular}
\end{tabular}
	\caption{Test case (b):
	Selected parameters and time nodes in PREIM for $\epsilon_{\textsc{eim}} = 10^{-1}$ (top) and $\epsilon_{\textsc{eim}} = 10^{-3}$ (bottom). The gray cells correspond to a new parameter selection and, therefore, 	to a new HF computation.}
	\label{b:mu_tb_1}
\end{table}
\begin{table}[htb]
	\centering
		\begin{tabular}{|c||c|c|c|c|c|c|c|c|c|c|c|c|c|c|}
		\hline
      $ m$ & $1$&$2$& $9$&$17$&$18$&$20$\\ \hline
      RB size & $5$ & $6$ &$6$& $7$& $9$ & $9$  \\ \hline
		\end{tabular}
	\caption{Test case (b): Size of RB generated within PREIM
	for $ \epsilon_\textsc{pod} = 5\cdot10^{-2}$; for $\epsilon_{\textsc{eim}} = 10^{-1}$, one stops at $m=9$, and for $\epsilon_{\textsc{eim}} = 10^{-3}$, one stops at $m=20$.}
	\label{b:pr_rb_sz}
\end{table}
\bl{We now investigate the performance of PREIM, which we run with thresholds 
$ \epsilon_\textsc{pod} = 5\cdot10^{-2}$ and either 
$\epsilon_\textsc{eim} = 10^{-1}$ or $\epsilon_\textsc{eim} = 10^{-3}$.
Table~\ref{b:mu_tb_1} shows the selected parameters and time nodes at each iteration. For $\epsilon_{\textsc{eim}} = 10^{-1}$, PREIM performs 9 iterations, and three parameters are selected for HF computations, whereas for $\epsilon_{\textsc{eim}} = 10^{-3}$, PREIM performs $11$ further iterations and six more HF computations to reach the requested threshold.
Moreover, the evolution of the size of the reduced basis within PREIM is shown in Table~\ref{b:pr_rb_sz}. As can be noticed, the approximation of the nonlinearity requires more computational effort than that of the solution manifold.}

\begin{figure}[htb]
\includegraphics[scale=0.65]{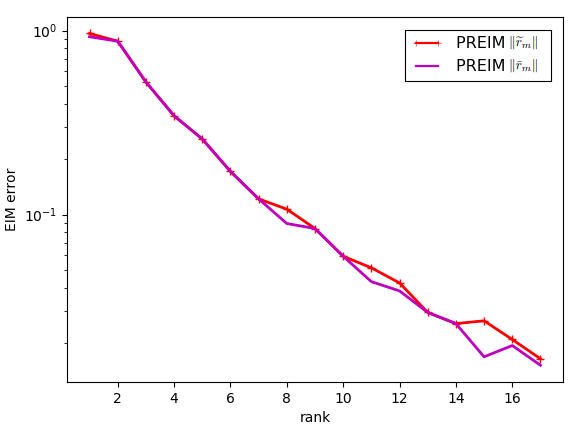}
\centering
\caption{Test case (b): EIM approximation errors $\|\widetilde r_m\|_{\ell^\infty(\Omega^{\mathrm{tr}})}$ and $\|\bar r_m\|_{\ell^\infty(\Omega^{\mathrm{tr}})}$ as a function of $m$ for PREIM with $\epsilon_\textsc{pod}=5\cdot10^{-2}$ and $\epsilon_\textsc{eim}=10^{-3}$.}
\label{b:r_rbar}
\end{figure}
\begin{figure}[htb]
\includegraphics[scale=0.65]{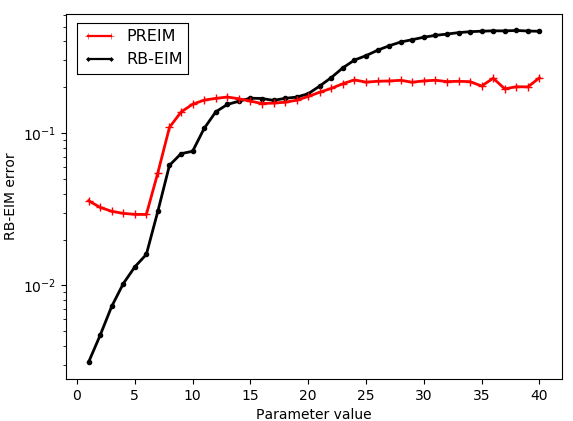}
\quad
\includegraphics[scale=0.64]{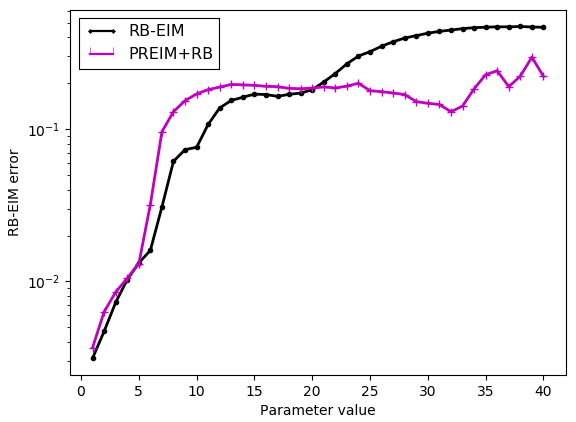}
\centering
\caption{Test case (b): RB approximation error $\|u_\mu -\hat u_\mu\|
_{\ell^2(I^\mathrm{tr};H^1(\Omega^\mathrm{tr}))}$ for 
$\epsilon_\textsc{pod}=5\cdot10^{-2}$ and $\epsilon_\textsc{eim}=10^{-3}$.}\label{b:combi}
\end{figure}
\begin{figure}[htb]
\includegraphics[scale=0.65]{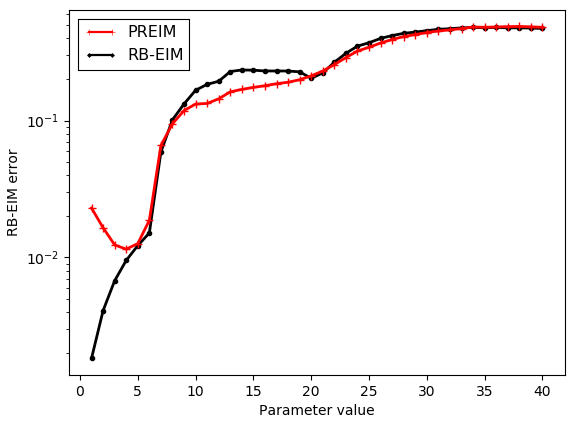}
\quad
\includegraphics[scale=0.64]{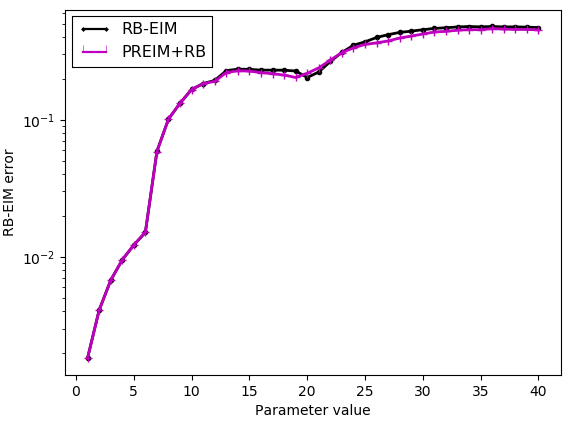}
\centering
\caption{Test case (b): RB approximation error $\|u_\mu -\hat u_\mu\|
_{\ell^2(I^\mathrm{tr};H^1(\Omega^\mathrm{tr}))}$ for 
$\epsilon_\textsc{pod}=2.5\cdot10^{-2}$ and $\epsilon_\textsc{eim}=10^{-4}$.}
\label{b:combi_conv}
\end{figure}
\bl{Figure~\ref{b:r_rbar} shows the decrease of the EIM approximation 
error on the nonlinearity for PREIM with $\epsilon_\textsc{pod}=5\cdot10^{-2}$ and $\epsilon_\textsc{eim}=10^{-3}$. We observe that each time a new HF trajectory is computed, i.e., whenever the quantities $\|\widetilde{r}_m\|_{\ell^\infty(\Omega^{\mathrm{tr}})}$ and $\|\bar r_m\|_{\ell^\infty(\Omega^{\mathrm{tr}})}$ differ, the difference is actually rather small, thereby confirming the already accurate approximation of the nonlinearity by the RB solutions in PREIM.
The left panel of Figure~\ref{b:combi} illustrates the space-time errors (measured in the 
$\ell^2(I^\mathrm{tr};H^1(\Omega^\mathrm{tr}))$-norm) on the trajectories produced by the standard RB-EIM and the PREIM procedures for the whole parameter range. 
We observe that for lower parameter values, PREIM delivers somewhat less accurate results, whereas the conclusion is reversed for higher parameter values. Altogether, both errors stay within comparable upper bounds. The right panel of Figure~\ref{b:combi} shows that the error component associated with the approximation of the nonlinearity is still the dominant one, except for the parameter range $[1,5]$ where the RB from the standard algorithm helps improve the error. Incidentally, we observe that these smaller values of the parameter were not selected within PREIM for approximating the nonlinearity. 
Finally, Figure~\ref{b:combi_conv} shows the same results for the tighter thresholds $\epsilon_{\textsc{pod}}=5\cdot10^{-2}$ and $\epsilon_{\textsc{eim}}=10^{-4}$.
Here, $14$ HF computations and $100$ PREIM iterations were needed. 
We can see that the PREIM error 
closely matches that of the standard RB-EIM procedure.}

\section{Conclusion and perspectives}\label{sec:conc}
In this work, we have devised a new methodology, called PREIM, 
that allows one to diminish the offline expenses incurred in the nonlinear RB method applied to unsteady nonlinear PDEs, 
\bl{as long as the computation of high-fidelity trajectories is the dominant part of the offline cost.} 
Numerical tests on two-dimensional nonlinear heat transfer problems with 
nonlinear thermal conductivities have \bl{illustrated} 
the computational efficiency 
and the accuracy of the algorithm.
In addition, the application of PREIM to an industrial test case of a three-dimensional
flow-regulation valve is ongoing.

\subsection*{Acknowledgments} This work is partially 
supported by Electricit\'e De France (EDF) and a 
CIFRE PhD fellowship from ANRT.
The authors are grateful to G. Blatman (EDF) and M. Abbas (EDF) for stimulating 
discussions and for their help on \texttt{code\_aster}. 
\bl{The authors are thankful to the two anonymous reviewers for
their valuable comments}.

\appendix

\section{Proper Orthogonal Decomposition}\label{pod_app}
The goal of this appendix is to briefly describe the procedure associated
with the notation 
\begin{equation} \label{app:pod}
(\theta_1,\ldots,\theta_N)=\mathrm{POD}(\mathcal S,\epsilon_\textsc{pod}),
\end{equation}
which is used in Algorithms~\ref{offline_rb}, \ref{progressive_pod}, and~\ref{preim},
where $\mathcal S = (v_1, \ldots, v_R)$ is composed of $R\ge1$ functions in the space $X$
and $\epsilon_\textsc{pod}$ is a user-prescribed tolerance. For simplicity,
we adopt an algebraic description\bl{, and we refer the reader to \cite{BeCOW:17} for further insight.} Let $(\varrho_1,\ldots,
\varrho_{\mathcal N})$ be a basis of $X$ where $\mbox{\rm dim} (X)
= \mathcal{N}$. For a function $w\in X$, we denote by 
$\mathbf w:=(w_j)_{1\leq j \leq \mathcal{N}}$
its coordinate vector in $\mathbb R^{\mathcal{N}}$, so that 
$w = \sum_{j=1}^{\mathcal{N}} w_j \varrho_j$. The algebraic counterpart
of~\eqref{app:pod} is that we are given
$R$ vectors forming the rectangular matrix 
$\mathbf S:=(\mathbf v_1, \ldots , \mathbf v_R)\in 
\mathbb R^{\mathcal{N}\times R}$, and we are looking for $N$ vectors
forming the rectangular matrix $\bm \Theta:=(\bm \theta_1, \ldots , \bm\theta_N)\in \mathbb R^{\mathcal{N}\times N}$. The vectors $(\bm \theta_1, \ldots , \bm\theta_N)$ are to be
orthonormal with respect to the Gram matrix of the inner product in $X$. 
In the present setting, we consider the Gram matrix 
$\mathbf C^\mathcal{N}\in \mathbb R^{\mathcal{N}\times \mathcal{N}}$ such that 
\begin{equation}
\mathbf C^\mathcal{N} = \Big(m(\varrho_n,\varrho_p)+ \eta a_0(\varrho_n,\varrho_p)\Big)_{1\leq p,n \leq \mathcal N },
\end{equation}
where $\eta>0$ is a user-prescribed weight and 
the bilinear forms $m$ and $a_0$ are defined in~\eqref{ma0}. Thus, we want to have
$\bm \theta_n^{\mathrm{T}} \mathbf C^\mathcal{N} \bm \theta_p=\delta_{n,p}$,
the Kronecker delta, for all $n,p\in\{1,\ldots,N\}$.

Let us set $\mathbf{T}:=( \mathbf C^\mathcal{N})^{\frac12} \mathbf S
\in \mathbb R^{\mathcal{N}\times R}$ 
and consider the integer $D = \min(\mathcal N, R)$ (\bl{in most situations}, we have $D=R$ and $D\ll \mathcal N$). 
Computing the Singular Value Decomposition~\cite{svd} of
the matrix $\mathbf T$, we obtain the real numbers
$\sigma_1 \geq \sigma_2 \geq \cdots \geq \sigma_D \geq 0$, the orthonormal family of
column vectors $(\bm{\xi}_n)_{1\leq n \leq D} \in (\mathbb R^{\mathcal{N}})^D$ 
(so that $\bm{\xi}_n^{\mathrm{T}}\bm{\xi}_p=\delta_{p,n}$) and 
the orthonormal family of column vectors  $(\hat{\bm{\psi}}_n)_{1\leq n \leq D} \in (\mathbb R^{R})^D$ 
(so that $\hat{\bm{\psi}}_n^{\mathrm{T}}\hat{\bm{\psi}}_p=\delta_{p,n}$), and we have
\begin{equation}\label{svd_eq}
\mathbf T = \sum_{n=1}^D \sigma_n \bm{\xi}_n \hat{\bm{\psi}}_n^{\mathrm{T}}.
\end{equation}
From~\eqref{svd_eq}, it follows that 
$\mathbf T \hat{\bm{\psi}}_n = \sigma_n \bm{\xi}_n$ and
$\mathbf T^{\mathrm{T}} \bm{\xi}_n = \sigma_n \hat{\bm{\psi}}_n$ for all
$n\in\{1,\ldots,D\}$. 
The vectors we are looking for are then given by 
$\bm{\theta}_n := ( \mathbf C^\mathcal{N})^{-\frac12}\bm{\xi}_n$
for all $n\in\{1,\ldots,N\}$ with 
$N:= \max \{ 1\leq n \leq D\ |\ \sigma_n \geq \epsilon_\textsc{POD} \}$.
It is well-known that the $N$-dimensional space spanned by the vectors
$(\bm\theta_n)_{1\le n\le N}$ minimizes the quantity
$\mathop{\inf}_{\mathbf z\in \mathbf Z_N} \sum_{r=1}^R (\mathbf v_r -\mathbf z)^{\mathrm{T}} \mathbf C^\mathcal{N} (\mathbf v_r- \mathbf z)$ 
among all the $N$-dimensional subspaces $\mathbf Z_N$ of $\mathbb R^{\mathcal{N}}$.
Moreover, we have
$\|v - \Pi_{Z_N}v\|_X\leq \sigma_{N+1}\|v\|_X$, 
for all $v \in \mathcal S$.

In practice, when $D=R$, we can avoid the computation of the matrix
$( \mathbf C^\mathcal{N})^{\frac12}$ and of its inverse by considering
the matrix of smaller dimension
$\mathbf T^{\mathrm{T}} \mathbf T = \mathbf S^{\mathrm{T}}\mathbf C^\mathcal{N} \mathbf S
\in \mathbb R^{R\times R}$. 
Solving for the eigenvalues of  $\mathbf T^{\mathrm{T}} \mathbf T$, 
we obtain the vectors $\hat{\bm{\psi}}_n$ with associated eigenvalues 
$\sigma_n^2$ since we have
$\mathbf T^{\mathrm{T}} \mathbf T \hat{\bm\psi}_n=\sigma_n\mathbf T^{\mathrm{T}} \bm\xi_n
 =\sigma_n^2\hat{\bm\psi}_n$.
Then, the vectors $(\bm\theta_n)_{1\leq n\leq N}$ are obtained as
$
\bm\theta_n = \left( \mathbf C^\mathcal{N}\right)^{-\frac12}\bm{\xi}_n
= \frac{1}{\sigma_n}\left( \mathbf C^\mathcal{N}\right)^{-\frac12}\mathbf T \hat{\bm{\psi}}_n 
= \frac{1}{\sigma_n} \mathbf S \hat{\bm{\psi}}_n.
$


\begin{thebibliography}{10}

\bibitem{eim04}
{\sc M.~Barrault, Y.~Maday, N.~C. Nguyen, and A.~T. Patera}, {\em An `empirical
  interpolation' method: application to efficient reduced-basis discretization
  of partial differential equations}, C. R. Math. Acad. Sci. Paris, 339 (2004),
  pp.~667--672.

\bibitem{BeCOW:17}
{\sc P.~Benner, A.~Cohen, M.~Ohlberger, and K.~Willcox}, eds., {\em Model
  Reduction and Approximation}, SIAM, Philadelphia, 2017.

\bibitem{ser15}
{\sc C.~Daversin and C.~Prud'homme}, {\em Simultaneous empirical interpolation
  and reduced basis method for non-linear problems}, C. R. Math. Acad. Sci.
  Paris, 353 (2015), pp.~1105--1109.

\bibitem{drohmann}
{\sc M.~Drohmann, B.~Haasdonk, and M.~Ohlberger}, {\em Reduced basis
  approximation for nonlinear parametrized evolution equations based on
  empirical operator interpolation}, SIAM J. Sci. Comput., 34 (2012),
  pp.~A937--A969.

\bibitem{aster}
{\sc {Electricit{\'e} de France}}, {\em Finite element {$\bf\it code\_aster$},
  analysis of structures and thermomechanics for studies and research}.
\newblock Open source on www.code-aster.org, 1989--2017.

\bibitem{ern_guermond}
{\sc A.~Ern and J.-L. Guermond}, {\em Theory and Practice of Finite Elements},
  Applied Mathematical Sciences, Springer New York, 2004.

\bibitem{grepl12}
{\sc M.~A. Grepl}, {\em Certified reduced basis methods for nonaffine linear
  time-varying and nonlinear parabolic partial differential equations}, Math.
  Models Methods Appl. Sci., 22 (2012), pp.~1150015, 40.

\bibitem{grepl07}
{\sc M.~A. Grepl, Y.~Maday, N.~C. Nguyen, and A.~T. Patera}, {\em Efficient
  reduced-basis treatment of nonaffine and nonlinear partial differential
  equations}, M2AN Math. Model. Numer. Anal., 41 (2007), pp.~575--605.

\bibitem{haasdonk13}
{\sc B.~Haasdonk}, {\em Convergence rates of the {POD}-greedy method}, ESAIM
  Math. Model. Numer. Anal., 47 (2013), pp.~859--873.

\bibitem{freefem}
{\sc F.~Hecht}, {\em FreeFem++, Third Edition, Version 3.0-1. User's Manual},
  LJLL, University Paris VI.

\bibitem{stamm16}
{\sc J.~S. Hesthaven, G.~Rozza, and B.~Stamm}, {\em Certified reduced basis
  methods for parametrized partial differential equations}, SpringerBriefs in
  Mathematics, Springer, Cham; BCAM Basque Center for Applied Mathematics,
  Bilbao, 2016.
\newblock BCAM SpringerBriefs.

\bibitem{Hinze_05}
{\sc M.~Hinze and S.~Volkwein}, {\em Proper orthogonal decomposition surrogate
  models for nonlinear dynamical systems: error estimates and suboptimal
  control}, in Dimension reduction of large-scale systems, vol.~45 of Lect.
  Notes Comput. Sci. Eng., Springer, Berlin, 2005, pp.~261--306.

\bibitem{podref}
{\sc K.~Kunisch and S.~Volkwein}, {\em Galerkin proper orthogonal decomposition
  methods for parabolic problems}, Numer. Math., 90 (2001), pp.~117--148.

\bibitem{MMOPR:00}
{\sc L.~Machiels, Y.~Maday, I.~B. Oliveira, A.~T. Patera, and D.~V. Rovas},
  {\em Output bounds for reduced-basis approximations of symmetric positive
  definite eigenvalue problems}, C. R. Acad. Sci. Paris S\'er. I Math., 331
  (2000), pp.~153--158.

\bibitem{eim09}
{\sc Y.~Maday, N.~C. Nguyen, A.~T. Patera, and G.~S.~H. Pau}, {\em A general
  multipurpose interpolation procedure: the magic points}, Commun. Pure Appl.
  Anal., 8 (2009), pp.~383--404.

\bibitem{svd}
{\sc B.~Noble and J.~W. Daniel}, {\em Applied Linear Algebra}, Prentice-Hall,
  3rd~ed., 1988.

\bibitem{reliable}
{\sc C.~Prud'Homme, D.~V. Rovas, K.~Veroy, L.~Machiels, Y.~Maday, A.~T. Patera,
  and G.~Turinici}, {\em {Reliable Real-Time Solution of Parametrized Partial
  Differential Equations: Reduced-Basis Output Bound Methods}}, {Journal of
  Fluids Engineering}, 124 (2001), pp.~70--80.

\bibitem{manzoni16}
{\sc A.~Quarteroni, A.~Manzoni, and F.~Negri}, {\em Reduced basis methods for
  partial differential equations}, La Matematica per il 3+2, Springer
  International Publishing, 2016.

\bibitem{rovas06}
{\sc D.~V. Rovas, L.~Machiels, and Y.~Maday}, {\em Reduced-basis output bound
  methods for parabolic problems}, IMA J. Numer. Anal., 26 (2006),
  pp.~423--445.

\bibitem{poim}
{\sc K.~Urban and B.~Wieland}, {\em Affine decompositions of parametric
  stochastic processes for application within reduced basis methods}, IFAC
  Proceedings Volumes, 45 (2012), pp.~716--721.

\end{thebibliography}
\end{document}